\documentclass[11pt]{article}

\usepackage{amsthm}
\usepackage{amsmath}
\usepackage{amssymb}
\def\ignore#1{\relax}

{\theoremstyle{plain}

}

{\theoremstyle{plain}

}

{\theoremstyle{plain}

}

{\theoremstyle{plain}

}

{\theoremstyle{plain}

}

{\theoremstyle{definition}

}

{\theoremstyle{definition}

}

{\theoremstyle{remark}

}

{\theoremstyle{remark}

}

\newcommand{\A}{ {\cal A} }

\newcommand{\C}{ {\bf C}  }

\newcommand{\FF}{ {\cal F} }

\newcommand{\V}{ {\cal V} }

\newcommand{\ncps}{ ( {\cal A} , \varphi ) }
\newcommand{\ncpsb}{ ( {\cal A} , \varphi , \V , f, \Phi ) }
\newcommand{\Kr}{{\rm Kr}}

\newcommand{\card}{ \mbox{card} }
\newcommand{\blockno}{ \mbox{blno} }
\newcommand{\Abs}{{\mbox{\rm Abs} }}

\newcommand{\id}{ id }
\newcommand{\inv}{ ^{-1} }

\newcommand{\ee}{ {\varepsilon}  }
\newcommand{\nca}{ NC^{(A)} }
\newcommand{\ncb}{ NC^{(B)} }

\newcommand{\ka}{ \kappa^{(A)} }
\newcommand{\kaa}{ \kappa^{(A')} }
\newcommand{\kb}{ \kappa^{(B)} }

\newcommand{\ThetaA}{ \Theta^{(A)} }
\newcommand{\ThetaB}{ \Theta^{(B)} }

\newcommand{\ecdef}{ \stackrel{def}{ \Leftrightarrow } }

\newcommand{\rstar}{ *_{r} }
\newcommand{\freestar}{ \framebox[7pt]{$\star$} }
\newcommand{\freestarA}{ \freestar^{(A)} }
\newcommand{\freestarB}{ \freestar^{(B)} }

\begin{document}

\title{\bf Non-crossing cumulants of type B}
\author{Philippe Biane \\
      Frederick Goodman \\
      Alexandru Nica }
\date{ }

\maketitle

\begin{abstract}
We establish connections between the lattices of non-crossing partitions
of type B introduced by V. Reiner, and the framework of the free probability
theory of D. Voiculescu.

Lattices of non-crossing partitions (of type A, up to now) have played
an important role in the combinatorics of free probability, primarily via
the non-crossing cumulants of R. Speicher. Here we introduce the concept of
{\em non-crossing cumulant of type B;} the inspiration for its definition
is found by looking at an operation of ``restricted convolution of
multiplicative functions'', studied in parallel for functions on symmetric
groups (in type A) and on hyperoctahedral groups (in type B).

The non-crossing cumulants of type B live in an appropriate framework
of ``non-commutative probability space of type B'', and are closely
related to a type B analogue for the R-transform of Voiculescu (which
is the free probabilistic counterpart of the Fourier transform). By
starting from a condition of ``vanishing of mixed cumulants of type B'',
we obtain an analogue of type B for the concept of free independence for
random variables in a non-commutative probability space.
\end{abstract}

\medskip

\medskip

\medskip

{\Large\bf Introduction}

\medskip

The free probability theory of D. Voiculescu is rooted in operator
algebras, but has connections to several other fields of mathematics,
and in  particular has a substantial combinatorial side.
The combinatorics of free probability is intimately related to the
Moebius inversion theory in the lattices of non-crossing partitions
first studied by G. Kreweras \cite{K}. The role of non-crossing
partitions in free probability was discovered by R. Speicher \cite{S1},
and was the subject of fairly intensive research after that. It is
noteworthy that the very concept of free independence for a family of
non-commutative random variables can be formulated in terms of the
concept of ``non-crossing cumulants'' introduced in \cite{S1}.

In another direction of development in the study of non-crossing
partitions, V. Reiner \cite{R} has introduced (with motivation from
problems on arrangements of hyperplanes) a type B analogue $\ncb (n)$
for the lattice $\nca (n)$ of non-crossing partitions of
$\{ 1, \ldots , n \}$.

$\ $

The content of the present paper can be summarized in one phrase as
follows: Starting from the lattices $\ncb (n)$ of Reiner, we work
out the framework of what should be a ``non-commutative probability
space of type B'', and, based on the type B analogue for non-crossing
cumulants, we propose a concept of free independence of type B.

$\ $

Let us now elaborate. As mentioned before, the concept of free independence
(of type A) can be described in terms of the non-crossing cumulants of
Speicher, via a condition called ``the vanishing of mixed cumulants''
(see review in Section 4.3 below). So if one wants to take a combinatorial
approach to free independence of type B, then the natural line of attack
goes by introducing the type B analogue of the non-crossing cumulants, and
by formulating the corresponding condition of vanishing of mixed cumulants.
But this line of attack presents us with an immediate difficulty, that
we don't have the framework where the non-crossing cumulants of type B
are to be defined. This framework (the ``non-commutative probability
space of type B'') has to be invented at the same time with the cumulants.

However, when looking at the theory in type A, one sees that the
non-crossing cumulants are closely related to the concept of R-transform
of Voiculescu (which is the free probabilistic counterpart of the
Fourier transform), and to a certain operation of ``boxed convolution'',
$\freestar$, on power series. So one can start the attack by defining the
type B analogue $\freestarB$ for $\freestar$. This turns out to be feasible
(in a canonical way, in fact) because of the following reason: there exists
a common idea, of studying a concept of ``convolution of multiplicative
functions'', which produces both $\freestar$ ( = $\freestarA$ ) and
$\freestarB$ by appropriate particularizations.
In the present paper we choose to present the convolution of multiplicative
functions in the framework of Cayley graphs of groups: the particularization
which leads to $\freestarA$ is obtained by considering symmetric groups, and
the one which leads to $\freestarB$ is obtained by considering
hyperoctahedral
groups. The approach in terms of Cayley graphs is made possible by the fact
that $\nca (n)$ and $\ncb (n)$ embed naturally into the Cayley graphs of
the symmetric group $S_n$, and respectively of the hyperoctahedral group
$W_n$. (The embedding of $\nca (n)$ was observed in \cite{B2}; for
its type B analogue, see Section 3.2 below.)

$\ $

When worked out explicitly, $\freestarB$ is an associative binary
operation on series with coefficients from $\C^2$. The property of
$\freestarB$ which is important for our investigation is described as
follows (see Section 5.3 below):

$\ $

Let ${\cal C}$ be the algebra structure on $\C^{2}$which is obtained
by identifying $( \alpha ' , \alpha '' ) \in \C^{2}$ with the
$2 \times 2$ matrix
$\left[ \begin{array}{cl} \alpha' & \alpha ''  \\
                          0     & \alpha '
\end{array} \right]$. Then $\freestarB$ can be viewed as a boxed
convolution of type A, but with coefficients in ${\cal C}$:
\[
\freestarB \ = \ \freestarA_{\, \cal C}. \hspace{3cm} (I)
\]

$\ $

The equality $(I)$ is the consequence of the following simple, but important
fact which takes place at the level of lattices of partitions: for every
$n \geq 1,$ the natural ``absolute value map''
$\mbox{Abs} : \ncb (n) \rightarrow \nca (n)$ is an $(n+1)$-to-1 cover (see
Section 1.4 below). As a result of this, the summations over $\ncb (n)$
which are involved in the definition of $\freestarB$ can be pushed forward
in a controlled way to summations over $\nca (n),$ and $(I)$ follows.

The fact that the algebra structure ${\cal C}$ on $\C^2$ plays a role in
considerations
about the lattices $\ncb (n)$ had already been noticed in Reiner's work
\cite{R}. Theorem 16 on page 217 of \cite{R} can, in fact, be viewed as
a result about $\freestarB$ (a characterization of $\freestarB$ on a part
of its domain of definition), which relates $\freestarB$ to the algebra
${\cal C}$. This result could also be used to obtain a derivation of $(I)$,
different from the one outlined in the preceding paragraph.

$\ $

The developments described above suggest that the non-crossing cumulants
of 
type B should be $\C^{2}$-valued, and should be defined by the formulas
which one would normally use for ${\cal C}$-valued cumulants of type A
(with ${\cal C}$ the algebra structure appearing in $(I)$). When this is
done, a suitable concept of non-commutative probability space of type B
arises at the same time with the definition of the cumulants of
type B. From this point on we can pursue the program suggested at the
beginning of the discussion: consider the condition of vanishing of mixed
cumulants of type B, and rephrase it in terms of moments, in order to arrive
to the definition of what is ``free independence of type B'' (cf. Section
7.2 below).

As a general comment, we note that the type B structures which come in
discussion seem to always be ``superimposed'' on their counterparts
of type A (rather than being some totally new objects).
This starts at the level of lattices, and goes all the way to
algebras of non-commuting variables, where a ``non-commutative probability
space of type B'' is essentially given by a representation  of a
non-commutative probability space of type A.

Concerning possible directions for future research:
The free probability of type A is a rich theory, and there are quite a
few of its aspects -- e.g. the free central limit theorem, the theory of
reduced free products, or the Fock space models for free independence --
for which it is certainly worth looking for type B analogues. It would also
be interesting if connections could be established between free probability
of type B and the line of research in non-commutative probability started by
Bozejko and Speicher \cite{BS} on ``$\varphi- \psi$ independence''.

$\ $

Following  the introduction, the paper is divided into seven sections.

In the first three sections we go over basic facts about the lattice
$\ncb (n)$ and present the embedding of $\ncb (n)$ into the Cayley graph
of the hyperoctahedral group $W_n$.

The Section 4 provides a brief review of some definitions and basic facts
which are commonly used in the combinatorics of free probability (of type
A),
and for which type B analogues will be developed in the Sections 5-7.

The Section 5 is devoted to the operations of boxed convolution $\freestarA$
and $\freestarB$ which are obtained by suitably particularizing the concept
of ``restricted convolution of multiplicative functions'' to symmetric and
respectively to hyperoctahedral groups.

The concept of non-crossing cumulant of type B is introduced in the
Section 6; in the same section we also present an equivalent, more
explicit description of this concept, and we point out what is the
type B analogue for the R-transform of Voiculescu.

Finally, in the Section 7 we study the condition of vanishing of mixed
cumulants of type B, and we arrive at the analogue of type B for the
concept of free independence.

The cross-referencing between sections is done by using the subsection
number (e.g. ``Proposition 3.3'' refers to the unique proposition stated
in the Section 3.3).

$\ $

$\ $

{\bf\large Acknowledgment.} The work presented here started from
discussions
between the three authors at MSRI Berkeley, during the special program in
Operator Algebras in 2000-2001. We gratefully acknowledge the hospitality
of MSRI during that period.

\newpage

\section{Non-crossing partitions of type A and of type B.}

\medskip


{\subsection {Review of non-crossing partitions of type A.}}

The partition lattice of a finite set $F$ has as elements the partitions
\setcounter{equation}{0}
of $F$ into disjoint non-empty subsets.  The non-empty subsets making up
a partition are called the blocks (or classes) of the partition.
If $p$ is a partition of $F$, and $a,b \in F$, we write $a \sim_p b$ to
denote that $a$ and $b$ are in the same block of $p$.
Partitions are ordered by reverse refinement:  $p \le q$ if $p$ is a
refinement of $q$, that is, if every block of $p$ is contained in a block
of $q$.  This partial order has a maximum element $1_F$, which has $F$
as its only block, and a minimum element $0_F$, in which every block is a
singleton. 

Now suppose that $F$ is totally ordered. A partition $p$ of $F$ is said to
be non-crossing if whenever $a < b < c < d$ in $F$, and $a \sim_p c$,
$b \sim_p d$, it follows that $b \sim_p c$. The set $\nca(F)$ of
non-crossing
partitions of $F$ is itself a lattice, when considered with the partial
order
induced from the partition lattice of $F$. The same $1_F$ and $0_F$ as
before
serve as maximal and respectively minimal element of $\nca (F).$ When $F$ is
the set
\begin{equation}
[n] \ := \ \{ 1,2, \ldots , n \} ,
\end{equation}
considered with the usual order, we write $\nca(n)$ instead of
$\nca ( \ [n] \ )$. The lattice $\nca(n)$ was first studied by Kreweras
\cite{K}.

One can give a recursive criterion for a partition to be non-crossing:
$p$ is a non-crossing partition of $F$ precisely when $p$ has a block $A$
which is an interval in $F$, and $p \setminus \{A\}$ is a non-crossing
partition of $F \setminus A$.

A geometric picture of the non-crossing condition is obtained by placing
the points of $F$ in order around  a circle.  Given a partition $p$ of
$F$, form for each block of $p$ the convex hull of the block (i.e. the
smallest convex set in the plane containing the points of the block).
The partition is non-crossing precisely when the convex hulls of
different blocks are non-intersecting.  This makes it clear that the
non-crossing condition is preserved under cyclic permutations of $F$.

$\nca(n)$ has a complementation map $\Kr$ (introduced by Kreweras) given
as follows. Consider the totally ordered set
$$J = \{1 < \bar 1 < 2 < \bar 2 < \cdots < n < \bar n\}.$$
For $p \in \nca(n)$, its complement $\Kr(p)$ is the largest element  $q$
of 
$$\nca(\{\bar 1, \bar 2, \dots, \bar n\})  \cong \nca(n)$$
such that $p \cup q$ is a non-crossing partition of $J$. Then
$\Kr$ is an order reversing bijection of $\nca(n)$. The following property
of $\Kr$ is also worth recording:
\begin{equation}
\blockno ( p ) \ + \ \blockno ( \Kr (p) ) \ = \ n+1, \ \
\forall \ p \in \nca (n),
\end{equation} 
where $\blockno (p)$ stands for the number of blocks of the partition $p$.

A left-hand version  $\Kr'$ of the Kreweras complement has the same
description as $\Kr$, but with $J$ replaced by
$$ J' = \{\bar 1 < 1 < \bar 2 < 2 \cdots < \bar n < n \}. $$
One has $\Kr' \circ \Kr = \id$ on $\nca(n)$.

Note that, via suitable identifications, one can talk about $\Kr$ and
$\Kr '$ on $\nca (F)$, where $F$ is any totally ordered set.

\medskip

\medskip

\subsection{ Non-crossing partitions of type B.}

The type $B$ analogue of the lattice of non-crossing partitions was
introduced by Reiner \cite{R}. Consider the totally ordered set:
\begin{equation}
[ \pm n ] \ :=  \ \{1 < 2 < \cdots < n < -1 < -2 < \cdots < -n \},
\end{equation}
with its inversion map $a \mapsto -a$.  One defines $ \ncb(n)$  to be the
subset of $\nca( \ [ \pm n ] \ ) \cong \nca(2n)$ consisting of partitions
which are invariant under the inversion map.

If $\pi \in \ncb (n)$, then the
blocks of $\pi$ are of two types: those which are inversion invariant,
and those which are not. From the non-crossing condition it is easily seen
that $\pi$ can actually have at most one block which is inversion invariant;
if this exists, it will be called the zero-block of $\pi$. The other
blocks of $\pi$ must come in pairs: if $X$ is a non-inversion invariant
block, then  $-X$ is also a block, different from $X$.

It is immediate that $\ncb(n)$ is a sublattice of $\nca( \ [ \pm n ] \ )$,
containing the minimal and maximal element of $\nca ( \ [ \pm n ] \ )$.

Furthermore, it is easily seen that $\ncb(n)$ is closed under the Kreweras
complements $\Kr$ and $\Kr'$ considered on $\nca( \ [ \pm n ] \ )$. When
restricted from $\nca ( \ [ \pm n ] \ )$ to $\ncb (n)$, the maps $\Kr$ and
$\Kr'$ will, therefore, yield two anti-isomorphisms of $\ncb (n)$, inverse
to each other, and which will also be called Kreweras complementation maps
(on $\ncb (n)$). Note that for $\pi \in \ncb (n)$ there is no ambiguity
about the meaning of ``$\Kr ( \pi )$'', no matter whether $\pi$ is viewed
as an element of $\ncb (n)$ or of $\nca ( \ [ \pm n ] \ )$.

Let us observe that Equation (1.2) gives us:
\begin{equation}
\blockno ( \pi ) \ + \ \blockno ( \Kr ( \pi ) ) \ = \ 2n+1, \ \
\forall \ \pi \in \ncb (n).
\end{equation} 
This has the following consequence: given $\pi \in \ncb (n)$, exactly one
of the two partitions $\pi$ and $\Kr ( \pi )$ has a zero-block. Indeed,
a partition in $\ncb (n)$ has a zero-block if and only if it has an odd
number of blocks; and the Equation (1.4) implies that exactly one of
$\pi$ and $\Kr ( \pi )$ has an odd number of blocks.

\medskip

\subsection{Absolute value of a non-crossing partition of type B.}

\noindent
{\bf Notation.} Let $n$ be a positive integer, and consider the sets
$[n]$ and $[ \pm n ]$ appearing in the preceding subsections (cf.
Eqns. (1.1) and (1.3)). We denote by $\Abs : [ \pm n ] \rightarrow [n]$
the absolute value map sending $\pm i$ to $i$, for $1 \leq i \leq n$.
Moreover, if $X$ is a subset of $[ \pm n ]$, we will use the notation
$\Abs (X)$ for the set $\{ \Abs (x) \ : \ x \in X \} \subset [n]$.

$\ $

\noindent
{\bf Proposition and Definition.}  {\em Let $n$ be a positive integer, and
let $\pi$ be a partition in $\ncb (n)$. Then the sets of the form:
\[
\Abs (X),  \mbox{ with $X$ a block of $\pi$}
\]
form a non-crossing partition $p$ of $[n]$. This $p \in \nca (n)$ will be
called the absolute value of $\pi$, and will be denoted $\Abs ( \pi )$.}

\begin{proof}
It is clear that $p = \Abs(\pi)$ is a partition of $[n]$.  In fact,
$a \sim_p  b \Leftrightarrow$ $(a \sim_\pi b  \hbox{ and } \break -a
\sim_\pi -b)$
\hbox{ or } $(a \sim_\pi -b  \hbox{ and }  -a \sim_\pi b)$  .

We have to check the non-crossing condition for $p$.  Suppose that
$a < b < c < d$ in $[n]$ and $a \sim_p  c$, $b \sim_p d$.
Then we have $a \sim_\pi c$ or $a \sim_\pi -c$ in $[\pm n]$, and
$b \sim_\pi d$ or $b \sim_\pi -d$ in $[\pm n]$.

Suppose $a \sim_\pi c$.  If $b \sim_\pi d$, then $b \sim_\pi c$ by the
non-crossing condition for $\pi$, so $b \sim_p c$.  If $b \sim_\pi -d$,
then, since $a < b < c < -d$ in $[\pm n]$, we have again $b \sim_\pi c$
by the non-crossing condition for $\pi$, so again $b \sim_p c$.

The discussion of the case when $a \sim_\pi -c$ is similar, and leads to
the same conclusion that $b \sim_p c$.
\end{proof}

\medskip

\medskip

\subsection{The absolute value is a (n+1)-to-1 cover.}

As counted by Kreweras \cite{K}, the number of partitions in $\nca (n)$ is
a Catalan number,
\[
\card ( \ \nca (n) \ ) \ = \ \frac{1}{n+1}
\left(  \begin{array}{c}  2n \\ n  \end{array} \right) .
\]
Reiner \cite{R} observes that in the type B case we have simply
\[
\card ( \ \ncb (n) \ ) \ = \
\left(  \begin{array}{c}  2n \\ n  \end{array} \right) ,
\]
and that, in fact, several formulas in type B are simpler than their
counterparts in type A, because of the absence of the factor of $1/(n+1)$.

In the present paper we will use the following fact, which gives a nice
interpretation for the relation between the cardinalities of $\nca (n)$ and
of $\ncb (n)$.

$\ $

\noindent
{\bf Theorem.} {\em Let $n$ be a positive integer. Then
$\pi \mapsto \Abs ( \pi )$ is an $(n+1)$-to-1 map from $\ncb (n)$ onto
$\nca (n)$.}

$\ $

We break the argument proving the theorem into several lemmas.

$\ $

\noindent
{\bf Lemma 1.} {\em Let $n$ be a positive integer. We have the relation:
\begin{equation}
\Kr ( \Abs ( \pi )) \ = \ \Abs ( \Kr ( \pi )), \ \ \pi \in \ncb (n).
\end{equation} }
[ Note: On the left-hand side of (1.5), ``$\Kr$'' denotes a
Kreweras complement in $\nca (n)$; while on the right-hand side  of (1.5),
``$\Kr$'' denotes a Kreweras complement in $\ncb (n)$. ]

\medskip

\noindent
{\bf Proof of Lemma 1.} Fix $\pi \in \ncb (n)$ about which we prove (1.5).

Let us observe that $\Abs ( \pi ) \cup \Abs ( \Kr ( \pi ))$ is a
non-crossing
partition of the set $1 < \overline{1} < \cdots < n < \overline{n}$. This
follows by applying  Proposition 1.3 to $\pi \cup \Kr ( \pi )$, which
is a non-crossing partition of the set
$$1 < \overline{1} < \cdots < n < \overline{n} < -1 < - \overline{1} <
\cdots < -n < - \overline{n} . $$
Since $\Kr ( \Abs ( \pi ))$ is maximal with the property that
$\Abs ( \pi ) \cup \Kr ( \Abs ( \pi ))$ is non-crossing, it follows that
we have the inequality:
\[
\Abs ( \Kr ( \pi )) \ \leq \ \Kr ( \Abs ( \pi )).
\]

In order to complete the proof of (1.5), it is then sufficient to check
that $\Abs ( \Kr ( \pi ))$ and $\Kr ( \Abs ( \pi ))$ have the same number
of blocks. From (1.2) we know that $\Kr ( \Abs ( \pi ))$ has
$n+1 - \mbox{blno}( \Abs ( \pi ))$ blocks. On the other hand, when we use
(1.4) and take into account that exactly one of $\pi$ and $\Kr ( \pi )$
has a zero-block, we obtain that the number of blocks of
$\Abs ( \Kr ( \pi ))$ is also equal to
$n+1 - \mbox{blno} ( \Abs ( \pi ))$.   \qed

$\ $

$\ $

\noindent
{\bf Lemma 2.} {\em Let $n$ be a positive integer. Suppose that
$X,Y,Z$ are non-empty subsets of $[ \pm n ]$, such that all of the
following hold:
\begin{enumerate}
\item
$Z= -Z$, $X \cap (-X) = \emptyset$, $Y \cap (-Y) = \emptyset$.
\item
The family of sets $Z,X$ and $-X$ is non-crossing.
\item
The family of sets $Z,Y$ and $-Y$ is non-crossing.
\item
$\Abs (X) =  \Abs (Y) \subset [n]$.
\end{enumerate}
Then either $X=Y$, or $X= -Y$.  }

\medskip

\noindent
{\bf Proof of Lemma 2.} Fix a $j \in \Abs(X) = \Abs (Y)$. By replacing
if necessary $X$ with $-X$ and $Y$ with $-Y$, we can assume without loss
of generality that $j \in X \cap Y$. The conclusion of the proof
then has to be that $X=Y$.

Draw $1<2< \cdots <n< -1< -2 < \cdots < -n$ around a circle, and cut out of
the circle the convex hull (boundary included) of the points belonging to
$Z$. Thus we cut out a convex $(2m)$-gon, where $\card (Z) = 2m$, with
$m \leq n$; and what remains of the circle is a union of $2m$ domains,
each of them bounded by a side of the $2m$-gon and by an arc of the circle.
Note that none of these $2m$ domains can contain a pair of points $i$ and
$-i$, with $1 \leq i \leq n$. (This is because connecting two points drawn
around the circle and belonging to the same domain cannot intersect the
convex $(2m)$-gon which was cut out; while the line connecting $i$ and
$-i$ does intersect the $(2m)$-gon.)

Now look at the domain (out of the $2m$ domains constructed in the
preceding paragraph) which contains the point $j$. Let $U \subset [ \pm n ]$
be the set of points drawn around the circle, and which belong to that
domain. From the hypothesis that $X$ and $Z$ don't cross we obtain that
$X \subset U$; similarly, the hypothesis that $Y$ and $Z$ don't cross gives
that $Y \subset U$. Finally, the observation made at the end of the
preceding paragraph shows that the absolute value function is injective on
$U$; therefore, the hypothesis $\Abs (X) = \Abs (Y)$ implies $X=Y$. \qed

$\ $

\noindent
{\bf Lemma 3.} {\em Let $n$ be a positive integer, and let $\pi, \rho$ be
in $\ncb (n)$. Suppose that:

\begin{enumerate}
\item
$\Abs ( \pi ) = \Abs ( \rho )$, and that

\item
$\pi$ and $\rho$ have a common zero-block $Z = -Z \subset [ \pm n ]$.
\end{enumerate}

Then $\pi = \rho$.}

\medskip

\noindent
{\bf Proof of Lemma 3.} Let $X$ be a non-zero block of $\pi$. Since
$\Abs ( \pi ) = \Abs ( \rho )$, there exists a block $Y$ of $\rho$ such
that $\Abs (Y) = \Abs (X)$. $Y$ is a non-zero-block of $\rho$ (because
$\rho$ has only one zero-block which is $Z$, with $\Abs (Z)$ disjoint
from $\Abs (X) = \Abs (Y)$ ). Thus Lemma 2 applies to $X,Y$ and $Z$,
and gives us that either $X=Y$ or $X=-Y$. In either case, the pair of
blocks $X$, $-X$ of $\pi$ must coincide with the blocks $Y$, $-Y$ of
$\rho$.
In this way we obtain that every block of $\pi$ also is a block of $\rho$,
and the conclusion $\pi = \rho$ follows.  \qed

$\ $

\noindent
{\bf Proof of the Theorem.} We will prove the inequality:
\begin{equation}
\card ( \ \{ \pi \in \ncb (n) \ : \ \Abs ( \pi ) = p \} \ ) \ \leq \
n+1, \ \forall \ p \in \nca (n).
\end{equation}
This will imply the statement of the theorem, because we know that
$\card ( \ \ncb (n) \ ) \ = \ (n+1) \cdot \card ( \ \nca (n) \ ).$

We fix $p \in \nca (n)$ about which we show that (1.6) holds. We denote
$\Kr (p) =: q.$ Let $A_{1}, \ldots , A_{k}$ be the list of the blocks of
$p$, and let $A_{k+1}, \ldots , A_{n+1}$ be the list of the blocks of $q$;
this notation can be used because we know that the total number of blocks
of $p$ and $q$ is $n+1$ ( cf. Eqn.(1.2) ).

Let $\pi \in \ncb (n)$ be such that $\Abs ( \pi ) = p$. Then we have
$\Abs ( \Kr ( \pi )) = q$, by  Lemma 1. Exactly one of $\pi$ and
$\Kr ( \pi )$ have a zero-block $Z = -Z \subset [ \pm n ]$. If $\pi$ has
a zero-block $Z$, then $\Abs (Z) = A_m$ for some
$1 \leq m \leq k$; while if $\Kr ( \pi )$ has a zero-block $Z$, then
$\Abs (Z) = A_m$ for some $k < m \leq n+1$. In either case,
we end by assigning to $\pi$ a number $m \in \{ 1, \ldots , n+1 \}$,
determined by the equality $\Abs (Z) = A_m$.

In this way we get a function
\begin{equation}
\Phi : \{ \pi \in \ncb (n) \ : \ \Abs ( \pi ) = p \} \ \rightarrow \
\{ 1, \ldots , n+1 \} ,
\end{equation}
defined by setting $\Phi ( \pi ) := m$, where $m$ is obtained from $\pi$
in the way described in the preceding paragraph.

But now, from  Lemma 3 it follows that the function $\Phi$ defined
above is injective. Indeed, let $\pi , \rho$ be in the domain of $\Phi$
(i.e.  they are partitions in $\ncb (n)$ such that
$\Abs ( \pi ) =  \Abs ( \rho ) = p$), and suppose that
$\Phi ( \pi ) = \Phi ( \rho ) = m \in \{ 1, \ldots , n+1 \}$. If
$m \leq k$, then Lemma 3 applies directly to give $\pi = \rho$.
If $m \geq k+1$, then Lemma 3 applies to give that
$\Kr ( \pi ) = \Kr ( \rho )$ -- but this still implies that $\pi = \rho$,
since $\Kr$ is one-to-one on $\ncb (n)$.

Finally, since the function $\Phi$ of (1.7) is injective,  its domain
can have at most $n+1$ elements; this is exactly (1.6).
\qed

$\ $

$\ $

\noindent
{\bf Remark.} The proof of the theorem actually tells us how to concretely
construct the partitions $\pi \in \Abs^{-1} (p) \subset \ncb (n),$ for a
given $p \in \nca (n):$ First we choose a block of either $p$ or $\Kr (p),$
which is to be ``lifted into a zero-block'' (of either $\pi$ or
$\Kr ( \pi )$ ). Then all the blocks of $\pi$ are completely determined
by the fact they must have blocks of $p$ as absolute values, and that
they cannot cross the chosen zero-block.

For instance, suppose that $p = \{ (1,2), \ (3,4) \} \in \nca (4).$
We have $\Kr (p) = \{ (1), \ (2,4),$
$(3) \}$. If we e.g. choose that
the block $(1,2)$ of $p$ is the one to be lifted into a zero-block, then
we get the partition $\pi = \{ (1,2,-1,-2), (3,4), (-3,-4) \} \in
\Abs^{-1}(p)$. While if we choose that the block $(1)$ of $\Kr (p)$ is
the one to be lifted into a zero-block, then we get a partition
$\pi \in \Abs^{-1} (p)$ which has $\Kr ( \pi )$ =
$\{ (1,-1), \ (2,4), \ (-2,-4), \ (3), \ (-3) \}$, and which must,
therefore,  be $\pi = \{ (1,-2), \ (-1,2), \ (3,4), \ (-3,-4) \} .$

\medskip

\medskip

\medskip

\section{The ``Cayley graph'' framework}

\medskip

\subsection{ Marked groups.}

We consider objects of the form $(G,T)$, where $G$ is a group and $T$ is a
\setcounter{equation}{0}
finite set of generators of $G$. Such a pair $(G,T)$ is sometimes called
``a marked group''. We will assume that $T$ does not contain the unit $e$
of $G$, and more importantly that:

(i) $T$ is closed under taking the inverse
$(x \in T \Rightarrow x^{-1} \in T)$; and

(ii) $T$ is closed under conjugation
$(x \in T, \ c \in G \Rightarrow c^{-1}xc \in T)$.

\medskip

\medskip

\subsection{ Word-length and distance on a marked group.}

Let $(G,T)$ be as in Section 2.1. For every element $e \neq a \in G$
we define its length $|a|$ as the smallest positive integer $n$
with the property that $a$ can be written as $a= x_{1} \cdots x_{n}$ with
$x_{1}, \ldots , x_{n} \in T$. By convention, the length of $e$ is
$|e| := 0$.

It is well-known and easy to prove that the length function
$| \cdot | : G \rightarrow  \{ 0,1,2,3, \ldots \ \}$
has the following properties:
\begin{equation}
\left\{    \begin{array}{lll}
|ab| \leq |a| + |b|, & \forall & a,b \in G   \\
                     &         &             \\
|a^{-1}| = |a| ,     & \forall & a \in G     \\
                     &         &             \\
|c^{-1}ac| = |a| ,   & \forall & a,c \in G.
\end{array}    \right.
\end{equation}

As a consequence, it is immediate that the formula:
\begin{equation}
d(a,b) \ := \ |a^{-1}b| \ = \ |ba^{-1}|
\ = \ |b^{-1}a| \ = \ |ab^{-1}|
\end{equation}
makes sense and defines a distance on $G$. Moreover, the distance $d$
is invariant under left and right translations on $G$ (i.e.
$d( c'ac'', c'bc'' ) = d(a,b)$ for all $a,b,c',c'' \in G$).

\medskip

\medskip

\subsection{ The partial order on a marked group }

Let $(G,T)$ be a marked group as in Section 2.1. We introduce a partial
order
on $G$ by declaring that for $a,b \in G$ we have
\begin{equation}
a \leq b \ \ecdef \ d(e,a) + d(a,b) = d(e,b)
\end{equation}
(where $e$ is the unit of $G$, and $d$ is the distance on $G$, as in
Section 2.2). It is immediately checked that the prescription (2.3) defines
indeed a partial order on $G$, which has the unit $e$ as (unique) minimal
element. It is also clear that for $a,b \in G$ we have the implication
\begin{equation}
a \leq b \ \Rightarrow \ |a| \leq |b|,
\end{equation}
the converse of which isn't generally true.

For any $a,b \in G$ such that $a \leq b$ we will use the natural interval
notation
\begin{equation}
[ a,b ] \ := \ \{ c \in G \ : \ a \leq c \leq b \} .
\end{equation}

When trying to understand the partial order on $G$, a useful concept is the
one of ``cover''. For $a,b \in G$ we say that $b$ covers $a$ if $a \leq b$
and if there are no elements of $G$ properly sitting between $a$ and $b$
(i.e. if the interval $[ a,b ]$ is reduced to just $a$ and $b$). As is
easily checked, an equivalent description for the fact that $b$ covers $a$
is ``$a \leq b$ and $|b| = |a|+1$''.

\medskip

\medskip

\subsection{Cayley graphs.}  \label{subsection Cayley graphs}

Let $(G,T)$ be a marked group as in Section 2.1.
The Cayley graph of $(G,T)$ is a graph whose vertices are the elements of
$G$, and whose edges are the two-element subsets $\{ a,b \} \subset G$ such
that $d(a,b) = 1$ (with $d$ defined as in Equation (2.2)).

The distance and the partial order on $G$, as discussed in the Sections 2.2
and 2.3, have natural interpretations in the Cayley graph of $(G,T).$ For
instance: for $a,b \in G$, the prescription used in (2.3) to define what
it means that $a \leq b$ can now be interpreted as saying that ``$a$ lies
on a geodesic from $e$ to $b$, in the Cayley graph of $(G,T)$''.

\medskip

\medskip

Cayley graphs are generally considered for marked groups $(G,T)$ with fewer
conditions imposed on $T$ than we had in the Section 2.1 (see e.g. Chapter
IV
of \cite{dlH}). For our purposes the framework of Section 2.1 is
nevertheless
appropriate, because we will only focus on the following two classes of
examples.

\medskip

\medskip

\subsection{Example: the symmetric groups.} \label{symmetric group as
marked group}

Let $n \geq 2$ be an integer, and let $S_{n}$ denote the symmetric group
on the set $[n]$ (i.e. the group of all permutations of
$[n] = \{ 1, \ldots , n \}$). The permutations $t \in S_{n}$ will be
usually written in cycle notation. ( E.g.
$t = (1,3,4)(2,6)(5) \in S_{6}$ is the permutation
$\left(  \begin{array}{cccccc}
1 & 2 & 3 & 4 & 5 & 6  \\ 3 & 6 & 4 & 1 & 5 & 2
\end{array} \right)$,
which partitions the set $\{ 1, \ldots , 6 \}$ into three orbits,
and was, therefore, written as a product of 3 cycles. One often omits the
cycles corresponding to orbits of cardinality 1, thus the same $t \in S_6$
may also appear written as $t= (1,3,4)(2,6).$ ) The term ``long cycle''
is used for a permutation $t \in S_{n}$ which has only one orbit,
necessarily equal to $[n]$.

Let $T_{n}$ be the set of all transpositions in $S_{n}$
(i.e. the set of permutations in $S_{n}$ which have one orbit of length
2 and $n-2$ orbits of length 1). It is clear that $(S_{n}, T_{n})$ is a
marked group, satisfying all the conditions discussed in Section 2.1.
The length function for $(S_{n},T_{n})$ is described by the formula:
\begin{equation}
|t| \ = \ n - ( \# \mbox{ of orbits of $t$} ), \ \ t \in S_{n},
\end{equation}
as is easily checked.

We remark that this is not the usual (Coxeter) generating set and length
function for the symmetric groups;  the usual generating set consists of
adjacent transpositions $(i, i+1)$ only, and is not invariant under
conjugation.  See, for example, \cite{H}.

Concerning the partial order on $S_n$ (defined as in Section 2.3), it is
worth pointing out how the concept of cover is explicitly described in
this example. It is easily checked that for $t_1, t_2 \in S_n$ we have:
\begin{equation}
t_2 \mbox{ covers } t_1 \ \Leftrightarrow \ \left\{
\begin{array}{l}
\mbox{$t_2 = t_1 r, \ $ where $r = (i,j) \in T_n$ is such that} \\
\mbox{$i$ and $j$ belong to different orbits of $t_1$}
\end{array}  \right.
\end{equation}
The effect of the right multiplication
with $r$ in the equality $t_2 = t_1 r$ of (2.7) is that the two orbits of
$t_1$ which contain $i$ and $j$ are united into one orbit of $t_2$ (which
thus contains both $i$ and $j$).

\medskip

\medskip

\subsection{Example: the hyperoctahedral groups.}
\label{Bn marked group}

Let $n$ be a positive integer, and let $W_{n}$ denote the
hyperoctahedral group with $2^{n}n!$ elements (or in other words, the
Weyl group of type $\mbox{B}_{n}$). The realization of $W_{n}$ which
we will
use in the present paper is as the group of permutations $\tau$ of the
set $[ \pm n ]$ (with $[ \pm n]$ as defined in Eqn.(1.3)), which have
the property that 
\begin{equation}
\tau ( -i ) \ = \ - \tau (i), \ \ 1 \leq i \leq n.
\end{equation}
Hence, we will view $W_{n}$ as a subgroup of $S_{\pm n}$, the symmetric
group on $[ \pm n ]$.

Every $\tau \in W_{n}$ decomposes as a product of cycles (since $\tau$ is in
particular an element of $S_{\pm n}$). Because of (2.8), we see that the
cycle decomposition of $\tau$ may contain two kinds of cycles: some which
are inversion invariant, and some which are not. The cycles of the
non-invariant kind must come in pairs (e.g. if $(1,2,-6,3)$ is a cycle of
$\tau$, then $(-1,-2,6,-3)$ must also be a cycle of $\tau$).

As generating set for $W_{n}$ we take the set of reflections $R_{n}$ which
consists of all transpositions $(i,-i)$, $1 \leq i \leq n,$ and of all the
products of two transpositions $(i,j)(-i,-j)$ where $i,j \in [ \pm n ]$
and $|i| \neq |j|.$ Then $( W_{n}, R_{n})$ is a marked group, and $R_n$ has
all the properties considered in Section 2.1. The length function for
$(W_{n}, R_{n})$ is described by the formula:
\begin{equation}
| \tau | \ = \ n -  \left(
\begin{array}{c} 
\mbox{ \# of pairs of orbits of $\tau $}  \\
\mbox{ which are non-invariant under inversion}
\end{array}   \right) , \ \ \tau \in W_{n}.
\end{equation}
Again, this is not the usual (Coxeter) generating set and length function
on the hyperoctahedral groups;  the usual generating set consists of the
adjacent reflections \break $(i, i+1)(-i, -i-1)$ (for $1 \le i \le n-1$)
and the transposition $(1, -1)$ only, and is not invariant under
conjugation. 

Continuing the analogy with the preceding example, let us now pass to the
partial order on $W_n$ (defined as in Section 2.3), and let us look at the
explicit description for the concept of cover with respect to this partial
order. The description is more complicated than what we had in the
Example 2.5, but the reader should have no difficulty to check that it is
done as follows:

$\ $

\noindent
{\bf Lemma.} 
{\em Let $\tau_1 , \tau_2$ be in $W_n.$ We have that $\tau_2$ covers
$\tau_1$ in the partial order coming from $( W_n , R_n )$ if and only if
$\rho := \tau_1^{-1} \tau_2$ is in $R_n$ and falls in one of the following
four situations:

(a) $\rho = (i,-i)$, where $i$ and $-i$ belong to different orbits of
$\tau_1$.

(b) $\rho = (i,j)(-i,-j)$ with $|i| \neq |j|$ and where $i$ and $-i$ belong
to the same orbit of $\tau_1$, but $j$ and $-j$ do not belong to the same
orbit of $\tau_1$.

(c) $\rho = (i,j)(-i,-j)$ with $|i| \neq |j|$ and where no two of
$i,j,-i,-j$ belong to the same orbit of $\tau_1$.

(d) $\rho = (i,j)(-i,-j)$ with $|i| \neq |j|$, where $i$ and $-j$ belong
to the same orbit of $\tau_1$, and this orbit is not invariant under
inversion
(hence does not contain $-i$ and $j$).  }

$\ $

In the situations (a), (b), (c) of the preceding lemma, the effect of the
right 
multiplication with $\rho$ in ``$ \tau_2 = \tau_1 \rho$'' is that some
distinct orbits of $\tau_1$ are united to form larger orbits of $\tau_2$.
The situation (d) is different; in this case, if
$X \subset [ \pm n ]$ denotes the orbit of $\tau_1$ which contains $i$ and
$-j$ then the right multiplication with $\rho$ has the effect of replacing
the orbits $X$ and $-X$ of $\tau_1$ $( X \neq -X)$ by two inversion
invariant
orbits $Y$ and $Z$ of $\tau_2$ such that $Y \cup Z = X \cup (-X),$
$i,-i \in Y$ and $j,-j \in Z$.

An immediate consequence of the lemma is the following.

$\ $

\noindent
{\bf Corollary.} 
{\em If $\tau_1 \leq \tau_2$ in $W_n$, then $\tau_2$ has at least as many
inversion invariant orbits as $\tau_1$.}

\begin{proof} Without loss of generality we may assume that $\tau_2$
covers $\tau_1$. Then $\rho := \tau_1^{-1} \tau_2$ must fall in one of the
four situations (a)--(d) described in the preceding lemma. Let $N_1$ and
$N_2$ denote the number of inversion invariant orbits of $\tau_1$ and of
$\tau_2$, respectively. By direct inspection we see that $N_2 = N_1$ in the
situations (b) and (c), $N_2 = N_1 + 1$ in the situation (a), and
$N_2 = N_1 + 2$ in the situation (d). Hence, the inequality
$N_1 \leq N_2$ always holds.
\end{proof}

\medskip

\medskip

\subsection{Restricted convolution.}\label{section restricted covolution}

The last ingredient of the ``Cayley graph framework'' which we want to
consider is a convolution operation for complex-valued functions defined on
the corresponding group.

$\ $

\noindent
{\bf Definition.} Let $(G,T)$ be a marked group as in Section 2.1,
and let $\FF (G, \C )$ denote the set of all complex valued functions on
$G$.
On $\FF (G, \C )$ we define an operation of
restricted convolution $\rstar$ via the following formula:
\begin{equation} \label{definition restricted convolution}
(u \rstar v) (a) \ = \ \sum_{\begin{array}{c}
{\scriptstyle b,c \in G,\  bc=a,} \\
{\scriptstyle  \ |b|+|c|=|a|}
\end{array} }   u(b)v(c).
\end{equation}

$\ $

\noindent
{\bf Remarks.}
\begin{enumerate}
\item
On the right-hand side of Equation (2.10) we are dealing with a finite sum,
because $T$ is assumed to be finite (and thus for every $k \geq 0$ there
are at most $( \card (T) )^{k}$ elements $a \in G$ such that $|a| = k$).
\item
Referring to the partial order in a marked group, we could also write the
Equation (2.10) in the form:
\begin{equation} \label{second definition restricted convolution}
(u \rstar v) (a) \ = \ \sum_{b \in [e,a]} \ u(b)v(b^{-1}a).
\end{equation}
\item
While the Equations (2.10), (2.11) are indeed reminiscent of the convolution
operation on $G$, the condition $|b|+|c|=|a|$ required in the sum on the
right-hand side of (2.10) changes the things quite  a bit in comparison to
the unrestricted convolution. For instance it is immediate that we have
$( u \rstar v )(e) \ = \ u(e)v(e),$ or that
$( u \rstar v )(x)$ = $u(e)v(x) + u(x)v(e)$ for $x \in T$.
\end{enumerate}

$\ $

We leave it as an exercise to the reader to make the straightforward
verifications proving the following proposition.

$\ $

\noindent
{\bf Proposition} {\em
Let $(G,T)$ be a marked group. Then $\FF (G, \C )$
is a unital complex algebra with the usual (pointwise) vector space
operations and with the restricted convolution $\rstar$ as multiplication.
The unit of $\FF (G, \C )$ is the characteristic function $\chi_{e}$ of the
unit of $G$ ($\chi_{e} (e) =1$ and $\chi_{e} (a) = 0$ for $a \neq e$ in
$G$). }

\medskip

\medskip

\medskip

\section{Non-crossing partitions and Cayley graphs.}

\medskip
\subsection{$\nca(n)$ and the symmetric group $S_n$.}
\label{$\nca(n)$ and the symmetric group $S_n$.}
Consider again the marked group $(S_n , T_n )$ of
Section \ref{symmetric group as marked group}. We will look at the interval
$[e,c]$ in the corresponding partial order on $S_n$, where $e$ is the unit
of $S_{n}$ and  $c$ is the ``forward" long cycle $c = ( 1,2, \ldots , n)$ on
$[n]$. It was observed by one of us in \cite{B2} that this interval provides
a ``group-theoretic incarnation'' of $\nca (n)$. More precisely, we have:

$\ $

\noindent
{\bf Theorem} {\em (see \cite{B2}, Theorem 1).
Define  a map $\iota$ from the partition lattice of $[n]$ to  $S_n$, as
\setcounter{equation}{0}
follows: 

-- First define $\iota$ on the set of subsets of $[n]$ by
$\iota(\{a_1, a_2, \dots, a_k\}) = (a_1, a_2, \dots, a_k)$, if
$a_1 < a_2 < \cdots < a_k $;

-- Extend $\iota$ to partitions of $[n]$  by defining $\iota$ of a
partition to be the product of $\iota(F)$, where $F$ runs over the blocks
of the partition. (The product is well defined since the cycles $\iota(F)$
are disjoint and, hence, mutually commuting.)

Then the restriction of $\iota$ to $\nca(n)$ has range equal to the interval
$[e,c] \subset S_n$, and $\iota$ is an order isomorphism between
$\nca (n)$ and the interval $[e, c]$ in $S_n$. }

$\ $

$\ $

\noindent
{\bf Remark.} The Kreweras complements have an interesting
interpretation in terms of the poset isomorphism $\iota : \nca (n)
\rightarrow [e,c]$. One has
\begin{equation}
\iota( \Kr (p)) = \iota(p)\inv c, \mbox{ and }
\iota( \Kr '(p)) =  c \iota(p)\inv , \forall \ p \in \nca (n).
\end{equation}

$\ $

$\ $

It actually turns out that the order structure of every subinterval
$[ e,b ] \subset [ e,c ]$ (with $b \leq c$) is closely related to the
lattices of non-crossing partitions. This comes as a consequence of some
basic properties of the partial order on $S_n$, which were put into
evidence in \cite{B1} and are reviewed in the following proposition.

$\ $

\noindent
{\bf Proposition} {\em (see Section 1.2 in \cite{B1}).

\begin{enumerate}

\item
Suppose that $a \in S_n$ has word-length $|a|=k$, and that
$a = t_1 \cdots t_k$ is an expression for and element $a$ as a product of
$k$ transpositions. Then for every $1 \leq j \leq k$, the two elements of
$[n]$ which are  transposed by $t_j$ belong to the same orbit of $a$.

\item
If $a \leq b$ in $S_n$, then every orbit of $a$ is contained in an orbit
of $b$ (and as a consequence, every fixed point of $b$ must also be a fixed
point of $a$).

\item
Let $b \neq e$ be a permutation in $S_n$, let $F_1 , \ldots , F_k \subset
[n]$
be the orbits of $b$ which contain more than one element, and let
$b = b_1 \cdots b_k$ be the corresponding factorization of $b$ into cycles
(for every $1 \leq j \leq k$, $b_j$ acts like $b$ on $F_j$ and acts like
$e$ on $[n] \setminus F_j$). If $a_1 \in [ e, b_1 ], \ldots , a_k \in
[ e, b_k ]$, then the permutations $a_1, \ldots , a_k$ commute with each
other, and $a := a_1 \cdots a_k$ is in $[ e,b ]$.

\item
Let $e \neq b \in S_n$ and its factorization $b = b_1 \cdots b_k$ be as
in the preceding statement (3). Then every permutation $a \in [ e,b ]$ can
be 
uniquely written in the form $a = a_1 \cdots a_k$ where
$a_1 \in [ e, b_1 ], \ldots , a_k \in [ e, b_k ].$
 
\end{enumerate}  }

$\ $

Suppose now that $e \neq b \leq c$ (where $c = ( 1,2, \ldots , n)$ as
before)
and that $b$ is factored as $b = b_1 \cdots b_k$ as in the statements (3)
and (4) of
the preceding proposition. Then the proposition gives us that we have a
canonical poset isomorphism
\begin{equation}
[ e,b ] \ \cong \ [ e,b_1 ] \times \cdots \times [ e, b_k ].
\end{equation}

Moreover, let us fix a $j$, $1 \leq j \leq k$, and let us denote by $F_j$
the unique orbit of $b_j$ which is not reduced to one point. Then the
interval $[ e, b_j ]$ consists of permutations $a \in S_n$ such that $a$
fixes every element of $[n] \setminus F_j$, and such that we have
$a \mid F_j \leq b_j \mid F_j$
$( =  \ b \mid F_j )$, where the latter
inequality is considered in the appropriate marked group
$(S_{F_j}, T_{F_j})$. Since $b \mid F_j$ is a long cycle on $F_j$, the
theorem presented above in this subsection gives us that
$[ e, b_j ] \cong  \nca ( F_j ).$

So in conclusion, for  $e \neq b \leq c$, we obtain that:
\begin{equation}
[ e,b ] \ 
\begin{array}[t]{ll}
\cong & \nca ( F_1 ) \times \cdots \times \nca ( F_k )  \\
     &                                                 \\
\cong & \nca ( l_1 ) \times \cdots \times \nca ( l_k )  \\
\end{array}
\end{equation}
where $F_1 , \ldots , F_k$ are the orbits of $b$ consisting of more than
one element, and where $l_j$ is the cardinality of $F_j$, $1 \leq j \leq k.$

Since (due to the theorem presented above) $[ e,b ]$ is naturally identified
to a subinterval of $\nca (n)$, the isomorphism (3.3) is closely related to
some of the ``canonical factorizations'' for intervals of $\nca (n)$ which
are studied in Section 3 of \cite{S1}.

\medskip

\medskip

\subsection{$\ncb(n)$ and the hyperoctahedral group $W_n$.}
\label{$\ncb(n)$ and the hyperoctahedral group $W_n$}

We now take on the analogues of type B for the facts presented in
Section 3.1. We will use the Cayley graph framework of Section 2,
particularized to the case of the marked group $(W_n , R_n )$ which
appears in Section 2.6. We will denote by $\omega$ the long inversion
invariant cycle
\begin{equation}
\omega \ := \ (1,2, \ldots , n,-1,-2, \ldots , -n) \in W_n ,
\end{equation}
and we denote the unit of $W_n$ by $\ee$. In this subsection we will look
at the interval $[ \ee , \omega ] \in W_n$, considered with respect to the
partial order coming from $(W_n , R_n )$.

Observe that, since $\ncb(n) \subset \nca( [ \pm n ] )$, the map
$\iota : \nca( [ \pm n ]) \rightarrow S_{\pm n}$ (defined as in  Theorem
3.1) gives us by restriction a map $\iota : \ncb(n) \rightarrow S_{\pm n}$.
In fact, it is clear from the definitions that
$\iota(\ncb(n)) \subset W_n \subset S_{\pm n}$. We will prove the following
theorem:

$\ $

$\ $

\noindent
{\bf Theorem.} 
{\em The map $\iota$ from $\ncb(n)$ into $W_n$ is an order isomorphism
of $\ncb(n)$ onto the interval $[ \ee , \omega ] \subset W_n$ (where the
partial order on $W_n$ is as described in the Sections 2.3, 2.6). }

$\ $

During the proof of the theorem it will occasionally be convenient to use
facts about the marked group $( S_{\pm n}, T_{\pm n} )$ (defined in the way
indicated in Section 2.5, where $T_{\pm n}$ denotes the set of all
transpositions in $S_{\pm n}$). For an element $\tau \in W_n$ the
word-length
of $\tau$ in $(S_{\pm n}, T_{\pm n})$ is generally different from the
word-length of $\tau$ in $(W_n , R_n )$; in order to distinguish between
the two word-lengths, we will denote them as $| \tau |_A$ and $| \tau |_B$,
respectively.

We will use two lemmas.

$\ $

\noindent
{\bf Lemma 1.}   {\em  For $\pi \in \ncb(n)$,
\begin{enumerate}

\item
$\iota(\pi) \iota(\Kr(\pi)) = \omega$,

\item
$|\iota(\pi)|_B + |\iota(\Kr(\pi))|_B = |\iota(\omega)|_B = n$, and

\item
$\iota (\pi ) \in
[ e , \omega ]$.

\end{enumerate} }

\noindent
{\bf Proof of Lemma 1.}

{\em 1.} Follows from Eqn.(3.1), with $\pi$ viewed as an element
of $\nca ( \pm n)$. ( Here we take into account that, as pointed out
in Section 1.2, $\Kr ( \pi )$ has the same meaning when $\pi$ is viewed as
an element of $\ncb (n)$ or as an element of $\nca ( \pm n )$. )

\vspace{6pt}

{\em 2.} Let $z(\pi) = 1$ if $\pi$ has a zero block, and $z(\pi) = 0$
otherwise; thus $z(\Kr(\pi)) = 1 - z(\pi)$. Let $f(\pi)$ denote the number
of pairs of non-inversion invariant blocks of $\pi$. We have:
\begin{equation*}
\begin{aligned}
|\iota(\pi)|_A &= 2n - z(\pi) - 2 f(\pi), \\
|\iota(\Kr(\pi))|_A &= 2n - 1 +  z(\pi) - 2 f(\Kr(\pi)),
\end{aligned}
\end{equation*}
so
$$
\begin{aligned}
2n -1 = |\omega|_A =  |\iota(\pi)|_A + |\iota(\Kr(\pi))|_A &= 4n -1 -2
(f(\pi) + f(\Kr(\pi)),\  \text{or}  \\
f(\pi) +
f(\Kr(\pi) &= n.
\end{aligned}
$$
Since $|\iota(\pi)|_B = n - f(\pi)$, we conclude
$$
|\iota(\pi)|_B + |\iota(\Kr(\pi))|_B = 2n - (f(\pi) +
f(\Kr(\pi)) = n = |\omega|_B.
$$

\vspace{6pt}

{\em 3.} This is a clear consequence of 1 and 2.  \qed

$\ $

$\ $

\noindent
{\bf Lemma 2.}  {\em
If $\tau_1 , \tau_2 \in W_n$, and $\tau_1 \le \tau_2 \le \omega$ with
respect to the partial order coming from $( W_n , R_n )$, then
$\tau_1 \le \tau_2$ with respect to the partial order coming from
$( S_{\pm n} , T_{\pm n} )$. }

\medskip

\noindent
{\bf Proof of Lemma 2.} Observe first that any
$\tau \in [ \ee , \omega ] \subset W_n$ has at most one inversion invariant
orbit. This follows from the Corollary 2.6, and the fact that $\omega$ has
one inversion invariant orbit.

We will prove the following implication (which clearly entails the statement
of the lemma): If $\tau_1 , \tau_2 \in \ [ \ee , \omega ]$ are such that
$\tau_2$ covers $\tau_1$ with respect to the order coming from
$( W_n , R_n )$, then we must have that $\tau_1 \leq \tau_2$ with respect to
the order coming from $( S_{ \pm n} , T_{ \pm n} )$.

So let us fix $\tau_1 , \tau_2 \in [ \ee , \omega ]$ such that $\tau_2$
covers $\tau_1$ with respect to the order coming from $(W_n , R_n )$.
Then $\rho := \tau_1^{-1} \tau_2$ must fall in one of the four situations
described in Lemma 2.6. Note that, in fact, $\rho$ cannot fall in the
situation (d) of that lemma; indeed (as pointed out in the discussion
following to Lemma 2.6), if $\rho$ would be in the situation (d) then it
would
follow that $\tau_2$ has at least two inversion invariant orbits, in
contradiction to the observation made at the beginning of this proof.

Hence, we have $\tau_2 = \tau_1 \rho$, with $\rho$ in one of the situations
(a), (b), or (c) described in Lemma 2.6. By comparing these situations
(a), (b), (c) with the equivalence stated in (2.7) of Section 2.5, it is
immediately seen that in all cases ((a), (b) and (c)) we will have indeed
that $\tau_1 \leq \tau_2$ with respect to the partial order coming from
$( S_{\pm n} , T_{\pm n} ).$   \qed

$\ $

$\ $

\noindent
{\bf Proof of the Theorem.} From Lemma 1 it follows that
$\iota(\ncb(n)) \subseteq [ \ee, \omega]$.

Next we have to show that if $\sigma \in [ \ee , \omega ]$, then there
exists $\pi \in \ncb(n)$ such that $\iota( \pi ) = \sigma$. This can be
done directly  by an inductive argument, but it is more convenient to
appeal to the corresponding result in type A. By  Lemma 2, we have
$\sigma \le \omega$ with respect to the order coming from
$(S_{\pm n} , T_{\pm n} )$, and, therefore, by Theorem 3.1, there exists
$\pi \in \nca( [ \pm n ] )$ such that $\iota( \pi ) = \sigma$. But since the
blocks of $\pi$ are the orbits of $\iota( \pi )$, we have $\pi \in \ncb(n)$.

We next pick partitions $\pi_1 \leq \pi_2$ in $\ncb (n)$, and we show
that $\iota ( \pi_1 ) \leq \iota ( \pi_2 )$ in $W_n$, with respect to the
order coming from $( W_n , R_n )$. Without loss of generality we can assume
that $\pi_2$ covers $\pi_1$ in $\ncb (n)$, i.e. that there are no elements
of $\ncb(n)$ properly between $\pi_1$ and $\pi_2$. A straightforward
inspection (also helped by  Proposition 2 in \cite{R}) shows that $\pi_1$
and $\pi_2$ must fall in one of the following three situations:

\vspace{6pt}

(a) $\pi_1$ has no inversion invariant block, and $\pi_2$ is obtained from
$\pi_1$ by merging a block and its inversion to form an inversion invariant
block. 

\vspace{6pt}

(b) $\pi_2$ is obtained from $\pi_1$ by merging the inversion invariant
block of $\pi_1$ with a pair of non-inversion invariant blocks.

\vspace{6pt}

(c) $\pi_2$ is obtained from $\pi_1$ by merging two non-inversion invariant
blocks, as well as merging the inversions of these two blocks.

\vspace{6pt}

\noindent
As is easily checked, the three situations listed above correspond exactly
(and in the same order) to the situations (a), (b), (c) described in
Lemma 2.6, and applied to the permutations $\iota ( \pi_1 )$ and
$\iota ( \pi_2 ).$ So in all the three situations we obtain that indeed
$\iota ( \pi_1 ) \leq \iota ( \pi_2 )$ (and, in fact even more, that
$\iota ( \pi_2 )$ covers $\iota ( \pi_1 )$ in the partial order coming from
$( W_n , R_n )$ ).

Finally, we have to show that if $\pi_1 , \pi_2 \in \ncb (n)$ and if
$\iota( \pi_1 ) \le \iota( \pi_2 )$ in $( W_n , R_n )$, then
$\pi_1 \leq \pi_2$ in $\ncb (n).$ Here again we can first invoke  Lemma 2
to obtain that $\iota( \pi_1 ) \le \iota( \pi_2 )$ in the partial order
coming from $( S_{\pm n} , T_{\pm n} )$, and then use  Theorem 3.1
to conclude that $\pi_1$ is a refinement of $\pi_2$.
\qed

$\ $

The bulk of this subsection was devoted to describing the interval
$[ \ee , \omega ] \subset W_n$, where
$\omega := ( 1,2, \ldots , n,-1,-2, \ldots , -n)$ is a long
cycle of invariant type. Let us conclude with a quick
look at the interval $[ \ee , \gamma ] \subset W_n$, where
$\gamma := (1,2, \ldots , n)(-1,-2, \ldots , -n) \in W_n$ is ``a long cycle
of non-invariant type''.

$\ $

$\ $

\noindent
{\bf Proposition.} {\em For $\gamma = (1,2, \ldots , n)(-1,-2, \ldots , -n)
\in W_n$, we have that $[ \ee , \gamma ] \cong \nca (n).$ A natural poset
isomorphism from $\nca (n)$ onto $[ \ee, \gamma ] \subset W_n$ is obtained
as follows: Let $\iota_o : \nca (n) \rightarrow \ncb (n)$ be the map
which associates to a partition $p = \{ F_1, \ldots , F_k \} \in \nca (n)$
the partition $\iota_o (p) := \{ F_1, \ldots , F_k, -F_1, \ldots , -F_k \}
\in \ncb (n).$ Let $\iota : \ncb (n)$
$\rightarrow W_n$ be the map defined at the beginning of
this subsection. Then the range of
$\iota \circ \iota_o : \nca (n) \rightarrow W_n$ is $[ \ee, \gamma ]
\subset W_n$, and $\iota \circ \iota_o$ induces a poset isomorphism between
$\nca (n)$ and $[ \ee , \gamma ]$. }

\medskip

\begin{proof}
Note that $\gamma \leq \omega$ (since $\omega = \gamma \cdot (n,-n)$, and
this is a covering situation as in (a) of Lemma 2.6). Hence,
$[ \ee , \gamma ] \subset [ \ee , \omega ],$ and the inverse of the poset
isomorphism $\iota : \ncb (n) \rightarrow [ \ee , \omega ]$ provided by
the preceding theorem will identify $[ \ee , \gamma ]$ with the interval:
\[
[ \ \{ (1),(2), \ldots ,(n),(-1),(-2), \ldots , (-n) \} \ , \
 \{ (1,2, \ldots ,n),(-1,-2, \ldots , -n) \} \ ] \subset \ncb (n).
\]
But the latter interval is precisely $\iota_o ( \nca (n) ).$
\end{proof}

\medskip

\medskip

\subsection{$\ncb(n)$ and the hyperoctahedral group $W_n$ (continued).}

In the same setting as in the Section 3.2, we will now establish the
type B analogue for the facts collected in  Proposition 3.1, and for the
isomorphism stated (in a type A setting) in the Equation (3.3) of the
Section 3.1.

In this section we will not need to maintain the separate notations for
the word-lengths $| \tau |_A$ and $| \tau |_B$ of an element $\tau \in W_n$
(as we did in Section 3.2); here ``$| \tau |$'' will always mean
``$| \tau |_B$''. We will also use the following notation:

$\ $

\noindent
{\bf Notation.} Let $\tau$ be in $W_n$, $\tau \neq \ee$. Let
\[
X_1, -X_1, \ldots , X_p, -X_p,Z_1, \ldots , Z_q
\]
be the list of distinct orbits of $\tau$ which have more than one element,
where $p,q \geq 0,$ and where $Z_j = -Z_j$ for $j \leq q.$ Denote:
\[
Y_1 = X_1 \cup (-X_1), \ldots , Y_p = X_p \cup ( -X_p), Y_{p+1} = Z_1,
\ldots , Y_{p+q} = Z_q,
\]
and for $1 \leq j \leq p+q$ let $\tau_j$ be the permutation in $W_n$
which acts like $\tau$ on $Y_j$ and like $\ee$ on $[ \pm n ] \setminus Y_j.$
Then the writing $\tau \ = \ \tau_1 \cdots \tau_{p+q}$
(with commuting factors $\tau_1 , \ldots , \tau_{p+q}$) will be called
``the cycle factorization of type B'' for $\tau$.

$\ $

In what follows we will focus on the situation when
$\tau \in [ \ee, \omega ],$ where $\omega$ is the long cycle appearing in
Section 3.2 (cf. Equation (3.4)). Note that, as a consequence of Corollary
2.6, a permutation $\tau$ in $[ \ee, \omega ]$ can have at most one
inversion 
invariant orbit; i.e. , for such a $\tau$ the parameter $q$ appearing in the
preceding notation is either 0 or 1.

$\ $

$\ $

\noindent
{\bf Proposition.} {\em

\begin{enumerate}

\item
Suppose that $\sigma \in [ \ee , \omega ] \subset W_n$ has word-length
$| \sigma |=k$, and that we have a writing
$\sigma = \rho_1 \cdots \rho_k$ with $\rho_1, \ldots , \rho_k \in R_n.$
Then for every $1 \leq j \leq k,$ any two elements of $[ \pm n ]$ which
are transposed by $\rho_j$ belong to the same orbit of $\sigma$.

\item
If $\sigma \leq \tau \leq \omega$ with respect to the partial order coming
from $( W_n , R_n )$, then every orbit of $\sigma$ is contained in an orbit
of $\tau$ (and as a consequence, every fixed point of $\tau$ must also be
a fixed point of $\sigma$).

\item
Let $\tau$ be in $[ \ee , \omega ]$, $\tau \neq \ee$, and let
$\tau = \tau_1 \cdots \tau_k$ be the cycle factorization of type B for
$\tau$.
If $\sigma_1 \in [ \ee, \tau_1 ], \ldots , \sigma_k \in [ \ee, \tau_k ]$,
then the permutations $\sigma_1, \ldots , \sigma_k$ commute with each other,
and $\sigma := \sigma_1 \cdots \sigma_k$ is in $[ \ee, \tau ]$.

\item
Let $\tau$ be in $[ \ee , \omega ]$, $\tau \neq \ee$, and let
$\tau = \tau_1 \cdots \tau_k$ be the cycle factorization of type B for
$\tau$.
Then every permutation $\sigma \in [ \ee, \tau ]$ can be uniquely written in
the form $\sigma = \sigma_1 \cdots \sigma_k$ where
$\sigma_1 \in [ \ee, \tau_1 ], \ldots , \sigma_k \in [ \ee, \tau_k ].$
 
\end{enumerate}  }

$\ $

In the proof of the proposition we will use the following simple lemma:

$\ $

$\ $

\noindent
{\bf Lemma.} {\em Let $Y_1, \ldots , Y_k$ be non-empty and pairwise disjoint
subsets of $[ \pm n ],$ such that $Y_j = -Y_j,$ $1 \leq j \leq k.$ Let
$\sigma_1 , \ldots , \sigma_k \in W_n$ be such that $\sigma_j$ fixes all the
elements of $[ \pm n ] \setminus Y_j$, $1 \leq j \leq k.$ Denote
$\sigma := \sigma_1 \cdots \sigma_k$ (commuting product). Then
$| \sigma | \ = \ | \sigma_1 |  + \cdots + | \sigma_k |.$ }

$\ $

\noindent
{\bf Proof of the Lemma.} This is immediate from the explicit formula for
$| \cdot |$ provided by Equation (2.9) of Section 2.6. \qed

$\ $

$\ $

\noindent
{\bf Proof of the Proposition.}

{\em 2.} Follows by combining  Lemma 2 of Section 3.2 with
part (2) of Proposition 3.1.

\vspace{10pt}

{\em 1.} Follows from  part (2) of the proposition, and the fact that for
every $1 \leq j \leq k$ we have $\rho_j \leq \sigma$. (The latter inequality
is easily proved directly from the definitions, by showing separately that
$\rho_j \leq \rho_1 \cdots \rho_j$ and that
$\rho_1 \cdots \rho_j \leq \sigma$.)

\vspace{10pt}

{\em 3.} For every $1 \leq j \leq k,$ let us denote by $Y_j = - Y_j$ the
subset of $[ \pm n ]$ where $\tau_j$ acts non-trivially. From  part (2)
of the proposition and the hypothesis that
$\sigma_j \leq \tau_j$ it follows that $\sigma_j$ fixes all the elements of
$[ \pm n ] \setminus Y_j$, $1 \leq j \leq k.$ Hence, the preceding lemma
can be applied to $\sigma_1 , \ldots , \sigma_k$, and gives us that:
\begin{equation}
| \sigma | \ = \ | \sigma_1 |  + \cdots + | \sigma_k |.
\end{equation}
Since it is also clear that $\tau_j$ and $\sigma_j^{-1} \tau_j$ fix all the
elements of $[ \pm n ] \setminus Y_j,$ $1 \leq j \leq k,$ the same lemma
also gives that:
\begin{equation}
| \tau | \ = \ | \tau_1 |  + \cdots + | \tau_k |
\end{equation}
and
\begin{equation}
| \sigma^{-1} \tau | \ = \ | \sigma_1^{-1} \tau_1 |  + \cdots +
| \sigma_k^{-1} \tau_k |.
\end{equation}

Now, the fact that $\sigma_j \in [ \ee , \tau_j ]$ is equivalent to the one
that $| \sigma_j | + | \sigma_j^{-1} \tau_j | = | \tau_j |,$
$1 \leq j \leq k.$ By adding together these equalities for
$1 \leq j \leq k,$ and by taking (3.5), (3.6), (3.7) into account, we obtain
that $| \sigma | + | \sigma^{-1} \tau | = | \tau |,$ i.e. that
$\sigma \in [ \ee , \tau ].$

\vspace{10pt}

{\em 4.} Let $Y_1, \ldots , Y_k \subset [ \pm n ]$ have the same
significance 
as in the proof of part (3) (for every $1 \leq j \leq k,$ $\tau_j$ acts like
$\tau$ on $Y_j$ and like $\ee$ on $[ \pm n ] \setminus Y_j$ ). Since every
orbit of $\sigma$ is contained in an orbit of $\tau$ (by  part (2) of the
proposition), it is immediate that $\sigma$ can be uniquely factored as
$\sigma = \sigma_1 \cdots \sigma_k$, where $\sigma_j$ fixes the elements of
$[ \pm n ] \setminus Y_j,$ $1 \leq j \leq k.$ Then the preceding lemma
applies,
and gives us equalities stated exactly as in (3.5), (3.6), (3.7) from the
proof of part (3). But then we have:
\[
| \tau |  \ = \ | \sigma | + | \sigma^{-1} \tau | \ \
\mbox{ (because $\sigma \in [ \ee , \tau ]$ )}
\]
\[
= \ \sum_{j=1}^{k} | \sigma_j | \ + \ \sum_{j=1}^{k} | \sigma_j^{-1} \tau_j
| 
\ \ \mbox{ (by (3.5), (3.7) )}
\]
\[
= \ \sum_{j=1}^{k} \Bigl( \ | \sigma_j | + | \sigma_j^{-1} \tau_j | \
\Bigr)
\]
\[
\geq \  \sum_{j=1}^{k} | \tau_j | \ \
\mbox{ (by the triangle inequality for $| \cdot |$)}
\]
\[
= \ | \tau | \ \ \mbox{ (by (3.6));}
\]
this can happen only if each of the inequalities
$| \sigma_j | + | \sigma_j^{-1} \tau_j | \geq | \tau_j |$ holds with
equality,
i.e. if $\sigma_j \in [ \ee , \tau_j ]$ for every
$1 \leq j \leq k.$ \qed

$\ $

\noindent
{\bf Remark.} Suppose now that $\tau \in [ \ee , \omega ],$ and
$\tau \neq \ee .$ Let $\tau = \tau_1 \cdots \tau_k$ be the cycle
factorization 
of type B for $\tau$, and for every $1 \leq j \leq k$ let $Y_j = -Y_j$
denote the subset of $[ \pm n ]$ where $\tau_j$ acts in a non-trivial 
way.  From  parts (3) and (4) of the preceding proposition, we get a 
canonical poset isomorphism:
\begin{equation}
[ \ee , \tau ] \ \cong \
[ \ee , \tau_1 ] \ \times \cdots \times [ \ee , \tau_k ].
\end{equation}
On the other hand, for every $1 \leq j \leq k$ it is immediate (by using
part (2) of the preceding proposition) that $[ \ee , \tau_j ]$ can be
identified 
with an interval going from the identity to a ``long cycle'' (either
inversion-invariant or not inversion-invariant) in the
hyperoctahedral group with symbols from $Y_j.$ By the results of Section
3.2,
the latter interval is in turn identified canonically to either
$\ncb ( \ \mbox{card}(Y_j)/2 \ )$ or $\nca ( \ \mbox{card}(Y_j)/2 \ )$
(depending on whether $Y_j$ was an orbit of $\tau$, or the union of two
disjoint orbits of $\tau$, inverse to each other in $[ \pm n]$ ).

The conclusion of this discussion is the following: If the given $\tau$ has
an inversion invariant orbit $Z = -Z,$ and if the other orbits with more
than 
one element for $\tau$ are denoted as $X_1, -X_1, \ldots , X_p, -X_p,$ then
\begin{equation}
[ \ee , \tau ] \ \cong \ \nca ( \ \mbox{card}( X_1 ) \ ) \times \cdots
\times \nca ( \ \mbox{card}( X_p ) \ ) \times
\ncb ( \ \mbox{card}( Z )/2 \ ).
\end{equation}
If we are in the opposite case (when $\tau$ has no inversion invariant
orbits), 
then we denote the orbits with more than one element for $\tau$ as
$X_1, -X_1, \ldots , X_k, -X_k,$ and we just get:
\begin{equation}
[ \ee , \tau ] \ \cong \ \nca ( \ \mbox{card}( X_1 ) \ ) \times \cdots
\times
\nca ( \ \mbox{card}( X_k ) \ ).
\end{equation}
The Equations (3.9) and (3.10) represent the type B analogue for the
Equation (3.3) of Section 3.1.

\medskip

\medskip

\medskip

\section{Review of basic definitions and facts from combinatorial free
probability (of type A).}

In this section we give a brief glossary of basic definitions and facts
which are used in the combinatorics of free probability, and for which
``type B analogues'' will be proposed in the following sections.

\medskip
\subsection{Non-commutative probability space.} The simplest concept of a
\setcounter{equation}{0}
``non-commutative probability space'', focusing only on algebraic and
combinatorial aspects, consists of a pair $\ncps$, where
$\A$ is a complex unital algebra (``the algebra of random variables''), and
where $\varphi : \A \rightarrow \C$ (``the expectation'') is a linear
functional, normalized by the condition that $\varphi (1) =1.$

\medskip

\medskip
\subsection{Free independence.} This concept is defined for a family of
unital subalgebras $\A_{1}, \ldots , \A_{k} \subset \A$, where
$\ncps$ is a non-commutative probability space. The precise definition
consists of a ``condition in moments'': $\A_{1}, \ldots , \A_{k}$ are freely
independent if and only if:
\begin{equation}
\left\{  \begin{array}{c}
\varphi ( a_{1}a_{2} \cdots a_{n} ) = 0   \\
                                        \\
\mbox{whenever $a_{1} \in \A_{i_{1}}, \ldots , a_{n} \in \A_{i_{n}}$} \\
                                        \\
\mbox{with $i_{1} \neq i_{2},i_{2} \neq i_{3}, \ldots ,i_{n-1} \neq i_{n}$}
\\
                                         \\
\mbox{and where $\varphi ( a_{1} ) = \cdots = \varphi (a_{n} ) = 0$}
\end{array} \right.
\end{equation}
(see e.g. Chapter 2 of \cite{VDN}).

\medskip

\medskip
\subsection{Non-crossing cumulants.} The non-crossing cumulant functionals
associated to a non-commutative probability space $\ncps$ were introduced
by R. Speicher \cite{S1}. They are a family of multilinear functionals
$$ \Bigl( \ka_{n} : \A^{n} \rightarrow \C \Bigr)_{n=1}^{\infty}$$
(where the superscript $A$ in ``$\ka_n$'' is a reminder that we are dealing
with objects ``of type A''). The equation which defines the functionals
$\ka_n$ is:
\begin{equation}
\sum_{p \in \nca(n)} \ \prod_{ F \ block \ of \ p} \
\ka_{\card (F)} \Bigl( \ (a_{1}, \ldots , a_{n}) \mid F \ \Bigr)
\ = \ \varphi (a_{1} \cdots a_{n}),
\end{equation}
holding for every $n \geq 1$ and for every $a_{1}, \ldots , a_{n} \in \A$.
In (4.2) we used the convention of notation that if
$F = \{ j_{1} < j_{2} < \cdots < j_{m} \}$ is a subset of
$\{ 1, \ldots , n \}$ and if $a_{1}, \ldots , a_{n} \in \A$, then
\begin{equation}
( a_{1}, \ldots , a_{n} ) \mid F \ := \
( a_{j_{1}}, a_{j_{2}}, \ldots , a_{j_{m}} ) \in \A^{m}.
\end{equation}

A recursive use of the Equation (4.2) gives explicit formulas for the
functionals $\ka_{n}$. For instance for $n=1,2,3$ we get:
\begin{equation}
\left\{  \begin{array}{lll}
\ka_{1} (a_{1}) & = & \varphi (a_{1}),  \\
                &   &                   \\
\ka_{2} (a_{1},a_{2}) & = & \varphi (a_{1}a_{2}) -
\varphi (a_{1}) \varphi (a_{2}),  \\
                &   &                   \\
\ka_{3} (a_{1},a_{2},a_{3}) & = &
{  \begin{array}[t]{cl}
\varphi (a_{1}a_{2}a_{3}) & - \varphi ( a_{1} ) \varphi ( a_{2}a_{3} ) -
\varphi ( a_{2} ) \varphi ( a_{1}a_{3} ) \\
                        & - \varphi ( a_{1}a_{2} ) \varphi (a_{3}) +
2 \varphi (a_{1}) \varphi (a_{2}) \varphi (a_{3}).
\end{array}  }
\end{array} \right.
\end{equation}

\medskip

The main reason for which non-crossing cumulants are an efficient tool in
free probability is that they provide a neat reformulation for the
definition
of free independence. More precisely, we have:

$\ $

\noindent
{\bf Proposition.} {\em Let $\ncps$ be a non-commutative probability space,
and let $\A_1, \ldots , \A_k$ be unital subalgebras of $\A$. Then
$\A_1, \ldots , \A_k$ are freely independent if and only if:
\begin{equation}
\left\{  \begin{array}{c}
\ka_{n} ( a_{1}, a_{2}, \ldots , a_{n} ) = 0   \\
                                        \\
\mbox{whenever $a_{1} \in \A_{i_{1}}, \ldots , a_{n} \in \A_{i_{n}}$} \\
                                        \\
\mbox{and $\exists \ 1 \leq s<t \leq n$ such that $i_{s} \neq i_{t}$.}
\end{array} \right.
\end{equation} }

$\ $

The condition in (4.5) is called ``vanishing of mixed cumulants''. The
equivalence between free independence and the vanishing of mixed cumulants
was first proved in \cite{S1}; an instructive alternative proof of this
equivalence appears in the Section 6 of \cite{NS3}.

\medskip

\medskip
\subsection{Moment series and R-transform.} Let $\ncps$ be a non-commutative
probability space, and let $a$ be an element of $\A$. The power series
\begin{equation}
M_a (z) := \sum_{n=1}^{\infty} \varphi (a^n) z^n
\end{equation}
is called the moment series of $a$ in $\ncps$, while the series
\begin{equation}
R_a (z) := \sum_{n=1}^{\infty} \ka_n (a, \ldots , a) z^n
\end{equation}
is called the R-transform of $a$ in $\ncps$. The concept of R-transform was
first studied by a method using Toeplitz matrices, in \cite{V1}, where the
following basic fact was proved: if $a,b \in \A$ and if the subalgebra of
$\A$ 
generated by $a$ is freely independent from the one generated by $b$, then
\begin{equation}
R_{a+b} (z) \ = \ R_a (z) + R_b (z).
\end{equation}
For the developments in the present paper it turns out to be more important
to look at the counterpart of (4.8) which expresses the R-transform of the
product $ab$ in terms of the the series $R_a$ and $R_b$ (under the same
hypothesis that the subalgebras generated by $a$ and $b$ are freely
independent). This multiplicative counterpart of (4.8) was analyzed in
\cite{NS2} (in fact, in the more general situation when one deals with
$k$-tuples $(a_1, \ldots , a_k)$ and $(b_1, \ldots , b_k)$ instead of just
$a$ and $b$) by introducing a certain operation of ``boxed convolution'',
which is the object of our next section.

\medskip

\medskip

\medskip

\section{ Boxed convolution of type A and of type B.}

\medskip
\subsection
{Review of boxed convolution of type A.}
\label{section type A boxed convolution}

\noindent
{\bf Definition}  \label{definition of free convolution of type A}
\begin{enumerate}

\item
We denote by $\ThetaA$ the set of power series
\setcounter{equation}{0}
of the form:
\begin{equation}
f( z ) \ = \ \sum_{n=1}^{\infty} \ \alpha_{n} z^{n} ,
\end{equation}
where the $\alpha_{n}$s are complex numbers.

\item 
On $\ThetaA$ we define a binary operation $\freestarA$, as follows. If
$f( z ) \ = \ \sum_{n=1}^{\infty} \ \alpha_{n} z^{n}$ and
$g( z ) \ = \ \sum_{n=1}^{\infty} \ \beta_{n} z^{n}$,
then $f \ \freestarA \ g$ is the series
$\sum_{n=1}^{\infty} \ \gamma_{n} z^{n}$, where
\begin{equation}
\gamma_{n}  \ := \  \sum_{ \begin{array}{c}
{\scriptstyle p \in \nca (n)} \\
{\scriptstyle p := \{ F_1, \ldots , F_k \} }  \\
{\scriptstyle Kr(p): = \{ E_1, \ldots , E_h \} }
\end{array}  } \ \
\Bigl( \prod_{j=1}^k \alpha_{\card(F_j)} \Bigr) \cdot
\Bigl( \prod_{i=1}^h \beta_{\card(E_i)} \Bigr) , \ n \geq 1.
\end{equation}

\end{enumerate}

\medskip

The operation $\freestarA$ was introduced in \cite{NS2}, in the more general
situation when we consider (instead of series as in (5.1)) series in $k$
non-commuting  indeterminants $z_1, \ldots ,z_k.$ For the sake of
simplicity,
we will limit the consideration of the present paper to the situation when
$k=1.$ The main point of $\cite{NS2}$ is that $\freestarA$ provides the
combinatorial description for the multiplication of two freely independent
elements, in terms of their R-transforms. More precisely, we have:

$\ $

\noindent
{\bf Theorem} {\em (see \cite{NS2}, Theorem 1.4). Let $\ncps$ be a
non-commutative probability space, and let $a,b \in \A$ be such that the
unital subalgebras of $\A$ generated by $a$ and by $b$ are freely
independent. Then the R-transform of the product $ab$ satisfies the
equation
\begin{equation}
R_{ab} \ = \ R_a \ \freestarA \ R_b.
\end{equation}  }

$\ $

The operation $\freestarA$ is associative, and has the identity series
$\Delta (z) =z$ as unit element. This can be checked either directly from
the combinatorial definition in (5.2) (as was done in \cite{NS2}), or by
exploiting the interpretation of $\freestarA$ provided by (5.3).

What we want to emphasize in our review here is that $\freestarA$ provides
the middle-ground between free probability (on one hand) and the Cayley
graph
framework of the Section 2 (on the other hand). The relation between
$\freestarA$ and free probability is illustrated by the preceding theorem.
In order to present the relation between $\freestarA$ and the Cayley graph
framework, we first record the remark that $\freestarA$ can be truncated to
an operation on $\C^n$, for every $n \geq 1.$

$\ $

\noindent
{\bf Remark.} Let 
$f( z ) \ = \ \sum_{n=1}^{\infty} \ \alpha_{n} z^{n}$ and
$g( z ) \ = \ \sum_{n=1}^{\infty} \ \beta_{n} z^{n}$ be two series in
$\ThetaA$. For every $n \geq 1,$ the coefficient $\gamma_n$ of order $n$
in the boxed convolution $f \ \freestarA \ g$ is defined by (5.2) as a
polynomial expression in $\alpha_1, \ldots , \alpha_n, \beta_1, \ldots ,
\beta_n.$ For instance for $n=1,2,3$ one gets:
\begin{equation}
\left\{  \begin{array}{l}
\gamma_1 \ = \ \alpha_1 \beta_1,   \\
                                  \\
\gamma_2 \ = \ \alpha_2 \beta_{1}^{2} + \alpha_{1}^{2} \beta_2, \\
                                  \\
\gamma_3 \ = \ \alpha_3 \beta_{1}^{3} + 3 \alpha_1 \alpha_2 \beta_1
\beta_2 + \alpha_{1}^{3} \beta_3.  \\
                                  \\
\end{array}  \right.
\end{equation}

As a consequence, for a fixed value of $n \geq 1,$ it makes sense to
define an operation $\freestarA_{\, n}$ on $\C^n$, which records how the
$n$-tuple of coefficients of order up to $n$ in $f \ \freestarA \ g$
is obtained from the corresponding $n$-tuples of coefficients in $f$
and in $g$. E.g. , the Equations (5.4) provide the explicit description
for the fact that
\[
( \gamma_1 , \gamma_2 , \gamma_3 ) \ = \
( \alpha_1 , \alpha_2 , \alpha_3 ) \ \freestarA_{\, 3} \
( \beta_1 , \beta_2 , \beta_3 ).
\]
\mbox{From} the properties of $\freestarA$ it follows that 
$\freestarA_{\, n}$ is associative and has unit equal to $(1,0, \ldots , 0)$
$\in \C^n$, for every $n \geq 1.$

$\ $

But now, the truncated operation $\freestarA_{\, n}$ turns out to be closely
related to the operation of restricted convolution in the Cayley graph
framework for the marked group $(S_n, T_n)$ discussed in Section 2.5.
This fact was observed in \cite{B1}, and is stated precisely as follows.

$\ $

\noindent
{\bf Proposition} {\em (see \cite{B1}, Section 3.1). Let $n$ be a fixed
positive integer. Consider the marked group $(S_n,T_n)$, as discussed in
Section 2.5. For every $n$-tuple $\alpha = ( \alpha_1, \ldots , \alpha_n)
\in \C^n$ we denote as $u_{\alpha} : S_n \rightarrow \C$ the function
defined by the following formula:
\begin{equation}
u_{\alpha} (t) \ = \ \alpha_1^{k_1(t)} \alpha_2^{k_2(t)} \cdots
\alpha_n^{k_n(t)}, \ \ t \in S_n ,
\end{equation}
where $k_{m}(t)$ stands for the number of orbits of cardinality $m$ of the
permutation $t$ ($1 \leq m \leq n$, $t \in S_n$).

Consider, on the other hand, the operation of restricted convolution
$\rstar$
for complex-valued functions on $S_n$ (as discussed in Section 2.7). Then
the set of functions $\{ u_{\alpha} \ : \ \alpha \in \C^n \}$ defined by
(5.5) is closed under $\rstar$, and the operation $\rstar$ on this set of
functions coincides with the truncated boxed convolution $\freestarA_{\,
n}$.
In other words: for every $\alpha , \beta \in \C^n$ we have that
$u_{\alpha} \ \rstar \ u_{\beta} = u_{\gamma},$ where
$\gamma = \alpha \ \freestarA_{\, n} \ \beta$.  }

$\ $

When put together, the preceding proposition and theorem give a connection
between 
free probability and the Cayley graph framework, which is obtained by using
the boxed convolution $\freestarA$ as an intermediate object.

We conclude this subsection by recording another fact about $\freestarA$
which will be useful in the sequel, namely that one can effortlessly define
``vector-valued versions'' of this operation.

$\ $

\noindent
{\bf Remark and Notation.} Let ${\cal C}$ be a unital commutative algebra
over $\C$. The formula (5.2) used in the definition of the operation
$\freestarA$ makes perfect sense if the ``scalars'' $\alpha_1, \alpha_2,
\alpha_3, \ldots$ and $\beta_1 , \beta_2, \beta_3 , \ldots$ appearing
there are elements of ${\cal C}$ (the sums and
products in (5.2) become sums and products in ${\cal C}$). When using
scalars from ${\cal C}$, one obtains a version of the boxed convolution of
type A which will be denoted as $\freestarA_{\, \cal C}$.

\medskip

\medskip

\subsection
{Boxed convolution of type B.}
\label{section type B boxed convolution}

We will now define the analogue of type B, $\freestarB$, of the operation
$\freestarA$ discussed in Section 5.1. In order to do so, we will repeat
the considerations which related $\freestarA$ to the Cayley graph
framework, in a context where we replace the symmetric groups by
hyperoctahedral groups. (The boxed convolution of type B will eventually
become the middle-ground between free probability of type B and the
Cayley graph framework; but this is not an issue for the moment, since
we haven't introduced the necessary elements of free probability of
type B.)

We start by looking at the analogues of type B for the functions
$u_{\alpha}$
which appeared in Equation (5.5) of Proposition 5.1. By taking into
consideration the specifics of the factorization into cycles in type B (as
discussed in the Notation 3.3) we come to the following:

$\ $

\noindent
{\bf Notation.} Let $n$ be a fixed positive integer, and consider the
hyperoctahedral group $W_n.$ Let $\alpha$ =
$( \ ( \alpha_1 ', \alpha_1 ''), \ldots , ( \alpha_n ', \alpha_n '') \ )$
be an $n$-tuple in $( \C^2 )^n.$ We denote as $u_{\alpha} : W_n \to \C$
the function defined by the formula
\begin{equation}
u_{\alpha} ( \tau ) \ = \
( \alpha_1 ')^{k_1( \tau )} ( \alpha_1 '')^{l_1( \tau )} \cdots
( \alpha_n ')^{k_n( \tau )} ( \alpha_n '')^{l_n( \tau )}, \ \
\tau \in W_n,
\end{equation}
where for every $\tau \in W_n$ and every $1 \leq m \leq n$: the number
$k_m( \tau )$ counts the pairs of orbits of $\tau$ which are not
inversion-invariant, and have cardinality $m$; and $l_m( \tau )$ counts the
orbits of $\tau$ which are inversion-invariant, and have cardinality $2m$.

$\ $

\noindent
{\bf Remark.} The next thing to do is examine the restricted convolution
of functions $u_{\alpha}$ of the kind introduced before. It is unfortunate
that a restricted convolution $u_{\alpha} \rstar u_{\beta},$ with
$\alpha , \beta \in ( \C^2 )^n,$ will not generally be a function
$u_{\gamma}$ for some $\gamma \in ( \C^2 )^n.$ This can be seen by direct
computation, already in the case when $n=2.$ Indeed, in the case when $n=2,$
let us denote the unit of $W_2$ as $\ee$, and let us denote
\[
\tau_1 = \left(  \begin{array}{rrrr}
 1  &  2  &  -1  &  -2  \\
-1  &  2  &   1  &  -2
\end{array}  \right) , \ \
\tau_2 = \left(  \begin{array}{rrrr}
 1  &  2  &  -1  &  -2  \\
 1  & -2  &  -1  &   2
\end{array}  \right) , \ \
\]
\begin{equation}
\tau \ = \ \tau_1 \tau_2 \ = \
 \left(  \begin{array}{rrrr}
 1  &  2  &  -1  &  -2  \\
-1  & -2  &   1  &   2
\end{array}  \right) .
\end{equation}
Then it is easily seen that every function $u_{\alpha}$ on $W_2$ must
satisfy the relation
\[
u_{\alpha} ( \ee ) u_{\alpha} ( \tau ) \ = \ u_{\alpha} ( \tau_1 )
\ u_{\alpha} ( \tau_2 ),
\]
while a convolution $u_{\alpha} \rstar u_{\beta}$ does not satisfy this
relation (unless some special conditions are imposed on $\alpha$ and
$\beta$). While concretely verifying the latter fact, the reader will
observe
that the source of the problem lies in the fact that the permutation $\tau$
in (5.7) can be factored not only as $\tau = \tau_{1} \tau_{2},$ but also as
\begin{equation}
\tau \ = \ \left(  \begin{array}{rrrr}
 1  &  2  &  -1  &  -2  \\
 2  &  1  &  -2  &  -1
\end{array}  \right)  \cdot
\left(  \begin{array}{rrrr}
 1  &  2  &  -1  &  -2  \\
-2  & -1  &   2  &   1
\end{array}  \right)
\end{equation}
(where both the factorizations in (5.7) and (5.8) use permutations from the
set of generators $R_2 \subset W_2$ considered in Section 2.6). Furthermore,
what makes the alternative factorizations in (5.7) and (5.8) coexist is
the fact that the permutation $\tau$ which is factored has more than one
inversion-invariant orbit.

Now, let us consider again the inversion-invariant long cycle $\omega$ which
played a prominent role in the Sections 3.2 and 3.3. We know that
any permutation  $\tau \in [ \ee , \omega ]$ can have at most one
inversion-invariant orbit (as consequence of the Corollary 2.6), so the
pathology described previously cannot occur for such $\tau$. And, in fact,
the 
following is true:

$\ $

\noindent
{\bf Proposition.} {\em Let $n$ be a positive integer, and let $\alpha$ =
$( \ ( \alpha_1', \alpha_1''), \ldots , ( \alpha_n', \alpha_n'') \ )$ and
$\beta = ( \ ( \beta_1', \beta_1''), \ldots , ( \beta_n', \beta_n'') \ )$
be in $( \C^2 )^n.$ Define the numbers $\gamma_1', \gamma_1'', \ldots ,
\gamma_n', \gamma_n'' \in \C$ by the formulas:
\begin{equation}
\gamma_m' \ = \ \sum_{ \begin{array}{c}
{\scriptstyle p \in \nca (m)} \\
{\scriptstyle p := \{ F_1, \ldots , F_k \} } \\
{\scriptstyle Kr(p) := \{ E_1, \ldots , E_h \} }
\end{array}  }  \ \
\Bigl( \prod_{j=1}^{k} \alpha_{\card(F_j)}' \Bigr) \cdot
\Bigl( \prod_{i=1}^{h} \beta_{\card(E_i)}' \Bigr) ,
\end{equation}
\begin{equation}
\gamma_m'' \ = \ \sum_{ \begin{array}{c}
{\scriptstyle \pi \in \ncb (m) \ with \ zero-block} \\
{\scriptstyle \pi := \{ Z, X_1, -X_1, \ldots , X_k, -X_k \} } \\
{\scriptstyle Kr( \pi ) := \{ Y_1, -Y_1, \ldots , Y_h, -Y_h \} }
\end{array}  }  \ \
\Bigl( \prod_{j=1}^{k} \alpha_{\card(X_j)}' \Bigr) \cdot
\alpha_{\card(Z)/2}'' \cdot
\Bigl( \prod_{i=1}^{h} \beta_{\card(Y_i)}' \Bigr)
\end{equation}
\[
+ \ \sum_{ \begin{array}{c}
{\scriptstyle \pi \in \ncb (m) \ without \ zero-block} \\
{\scriptstyle \pi := \{ X_1, -X_1, \ldots , X_k, -X_k \} } \\
{\scriptstyle Kr( \pi ) := \{ Z, Y_1, -Y_1, \ldots , Y_h, -Y_h \} }
\end{array}  }  \ \
\Bigl( \prod_{j=1}^{k} \alpha_{\card(X_j)}' \Bigr) \cdot
\beta_{\card(Z)/2}'' \cdot
\Bigl( \prod_{i=1}^{h} \beta_{\card(Y_i)}' \Bigr) ,
\]
for $1 \leq m \leq n$ (and where in (5.10) $Z$ stands for the zero-block
of either $\pi$ or $\Kr ( \pi )$ ). Then
\begin{equation}
( u_{\alpha} \rstar u_{\beta} ) \mid [ \ee , \omega ] \ = \
u_{\gamma} \mid [ \ee , \omega ] ,
\end{equation}
where $\gamma := ( \ ( \gamma_1', \gamma_1''), \ldots , ( \gamma_n',
\gamma_n'') \ ) \in ( \C^2 )^n$.  }

$\ $

In the proof of the proposition we will use the following

$\ $

\noindent
{\bf Lemma.} {\em Suppose that $Y_1, \ldots , Y_k$ are non-empty
inversion-invariant subsets of $[ \pm n ]$, such that
$Y_i \cap Y_j = \emptyset$ for $i \neq j,$ and suppose that
$\sigma_1, \ldots , \sigma_k \in W_n$ are such that $\sigma_j$ fixes all
the elements of $[ \pm n ] \setminus Y_j,$ $1 \leq j \leq k.$ Let $\alpha$
be in $( \C^2 )^n,$ and consider the function $u_{\alpha} : W_n \to \C$
defined as previously (cf. Equation (5.6)). Then
\begin{equation}
u_{\alpha} ( \sigma_1 ) \cdots u_{\alpha} ( \sigma_k ) \ = \
u_{\alpha} ( \sigma_1 \cdots \sigma_k ) \cdot u_{\alpha} ( \ee )^{k-1},
\end{equation}
where $\ee$ is the unit of $W_n.$ }

$\ $

The proof of the lemma is a straightforward application of  Equation (5.6)
defining $u_{\alpha},$ and is left to the reader.

$\ $

\noindent
{\bf Proof of the Proposition.} Note that in the case when $m=1,$ the
Equation (5.9) gives us that $\gamma_1' = \alpha_1' \cdot \beta_1'.$ This
shows in particular that $( u_{\alpha} \rstar u_{\beta} ) ( \ee )$ =
$u_{\gamma} ( \ee ),$ since
\[
( u_{\alpha} \rstar u_{\beta} ) ( \ee ) \ = \
u_{\alpha} ( \ee ) \cdot u_{\beta} ( \ee ) \ = \
( \alpha_1' )^n \cdot ( \beta_1' )^n,
\]
while $u_{\gamma} ( \ee ) = (\gamma_1')^n.$

For the rest of the proof we fix a permutation $\tau \in [ \ee , \omega ],$
$\tau \neq \ee$, about which we will prove that
$( u_{\alpha} \rstar u_{\beta} ) ( \tau )$ =
$u_{\gamma} ( \tau ).$

It is immediate that the value of $( u_{\alpha} \rstar u_{\beta} ) ( \tau )$
varies continuously as a function of
$( \alpha , \beta ) \in ( \C^2 )^n \times ( \C^2 )^n$; and also that
$u_{\gamma} ( \tau )$ depends continuously on $( \alpha , \beta )$ (since
$\gamma$ itself does so, as is clear from the Equations (5.9), (5.10) ).
By making if necessary a small perturbation of $\alpha_1'$ and $\beta_1'$
we can, therefore, assume, without loss of generality, that
$\alpha_1' \neq 0 \neq \beta_1'.$ Note that this implies
$u_{\alpha} ( \ee ) = {\alpha_1'}^n \neq 0,$
$u_{\beta} ( \ee ) = {\beta_1'}^n \neq 0.$

The permutation $\tau$ which we fixed has at most one inversion-invariant
orbit. Thus there are two cases to consider: when $\tau$ has no
inversion-invariant orbit, and when $\tau$ has exactly one
inversion-invariant
orbit $Z$. We will discuss the latter situation (the case without
inversion-invariant orbit is analogous, one just has to ignore the part
corresponding to $Z$ throughout the computations).

Let then $Z, X_1, -X_1, \ldots , X_k, -X_k$ be the list of orbits of $\tau$
which have more than one element. Let $\tau$ =
$\tau_0 \tau_1 \cdots \tau_k$ be the cycle decomposition of type B for
$\tau$
(as discussed in the Notation 3.3), where $\tau_0$ acts non-trivially on $Z$
and $\tau_j$ acts non-trivially on $X_j \cup (-X_j),$ $1 \leq j \leq k.$
With 
these notations it is clear that
\begin{equation}
u_{\gamma} ( \tau ) \ = \ \gamma_{\card(Z)/2} '' \cdot
\prod_{j=1}^k  \gamma_{\card(X_j)} ' \cdot
( \gamma_1 ')^{ n- card(Z)/2 - \sum_{j=1}^k card(X_j) },
\end{equation}
and what we have to prove is that $( u_{\alpha} \rstar u_{\beta} ) ( \tau )$
is also equal to the same quantity.

We compute:
\[
( u_{\alpha} \rstar u_{\beta} ) ( \tau ) \ = \
\sum_{\sigma \in [ \ee , \tau ]}
u_{\alpha} ( \sigma ) u_{\beta} ( \sigma^{-1} \tau )
\mbox{  (by the definition of $\rstar$) }
\]
\[
= \ \sum_{  \begin{array}{c}
{\scriptstyle \sigma_0 \in [ \ee , \tau_0], \ \ldots}   \\
{\scriptstyle \ldots , \sigma_k \in [ \ee , \tau_k]  }
\end{array}  } \ \ u_{\alpha} ( \sigma_0 \cdots \sigma_k )
u_{\beta} ( \sigma_0^{-1} \tau_0 \cdots \sigma_k^{-1} \tau_k )
\mbox{  (by Prop. 3.3) }
\]
\[
= \ \sum_{  \begin{array}{c}
{\scriptstyle \sigma_0 \in [ \ee , \tau_0], \ \ldots}   \\
{\scriptstyle \ldots , \sigma_k \in [ \ee , \tau_k]  }
\end{array}  } \ \ \Bigl( \prod_{j=0}^k u_{\alpha} ( \sigma_j ) \Bigr)
\cdot u_{\alpha} ( \ee )^{-k} \cdot
\Bigl( \prod_{j=0}^k u_{\beta} ( \sigma_j^{-1} \tau_j ) \Bigr)
\cdot u_{\beta} ( \ee )^{-k} \cdot
\]
(by the lemma preceding this proposition)
\begin{equation}
= \ \Bigl( u_{\alpha} ( \ee ) u_{\beta} ( \ee ) \Bigr)^{-k} \cdot
\prod_{j=0}^{k} \Bigl( \sum_{\sigma_j \in [ \ee , \tau_j ] } \
u_{\alpha} (\sigma_j) u_{\beta} ( \sigma_j^{-1} \tau_j ) \ \Bigr) .
\end{equation}
At this moment we pick a value of $j \in \{ 0,1, \ldots , k \}$ and we look
at the sum over $[ \ee , \tau_j ]$ which appeared in (5.14). Exactly as
explained at the end of Section 3.3, this sum is converted into a sum over
$\ncb ( \ card(Z)/2 \ ),$ if $j=0,$ and into a sum over
$\nca ( \ card(X_j) \ ),$ if $1 \leq j \leq k.$ It is immediately verified
that the change of variable  from
``$\sum_{\sigma_0 \in [ \ee , \tau_0 ] }$'' to
``$\sum_{\pi \in \ncb ( card(Z)/2 ) }$'' takes us to a sum as described
in Equation (5.10), thus leading to:
\begin{equation}
\sum_{\sigma_0 \in [ \ee , \tau_0 ] } \
u_{\alpha} (\sigma_0) u_{\beta} ( \sigma_0^{-1} \tau_0 ) \ = \
( \alpha_1 ' \beta_1 ' )^{n - card(Z)/2} \cdot \gamma_{\card(Z)/2} ''.
\end{equation}
(The power of $\alpha_1 ' \beta_1 '$ appearing on the right-hand side of
(5.15) comes from the $n - card(Z)/2$ pairs of fixed points which both
$\sigma_0$ and $\sigma_0^{-1} \tau_0$ have in $[ \pm n ] \setminus Z.$ )
Similarly, for $1 \leq j \leq k$ we get that:
\begin{equation}
\sum_{\sigma_j \in [ \ee , \tau_j ] } \
u_{\alpha} (\sigma_j) u_{\beta} ( \sigma_j^{-1} \tau_j ) \ = \
( \alpha_1 ' \beta_1 ' )^{n - card(X_j)} \cdot \gamma_{\card(X_j)} '.
\end{equation}

By replacing (5.15) and (5.16) in (5.14) we come to:
\[
( \ u_{\alpha} \rstar u_{\beta} \ ) ( \tau ) \ = \
( \ u_{\alpha} ( \ee ) u_{\beta} ( \ee ) \ )^{-k} \times
\]
\[
\times ( \alpha_1 ' \beta_1 ')^{n(k+1) - card(Z)/2 - \sum_{j=1}^k card(X_j)
}
\cdot \gamma_{\card(Z)/2} '' \cdot \prod_{j=1}^k \gamma_{\card(X_j)} '.
\]
This is indeed equal to the expression from (5.13), since
$( \ u_{\alpha} ( \ee ) u_{\beta} ( \ee ) \ )^{-k}$ =
$( \alpha_1 ' \beta_1 ')^{-nk}$, and $\alpha_1 ' \beta_1 ' = \gamma_1 '$.
\qed

$\ $

$\ $

We are thus led to introduce the type B analogue for the operation of boxed
convolution, in the following way.

$\ $

\noindent
{\bf Definition}  \label{definition of free convolution of type B}
\begin{enumerate}

\item
We denote by $\ThetaB$ the set of power series of the form:
\begin{equation}
f( z ) \ = \ \sum_{n=1}^{\infty} \ ( \alpha_{n}', \alpha_{n}'') z^{n} ,
\end{equation}
where the $\alpha_{n}'$s and $\alpha_{n}''$s are complex numbers.

\item 
Let $f(z) := \sum_{n=1}^{\infty}  ( \alpha_n ' , \alpha_n '') z^n$ and
$g(z) := \sum_{n=1}^{\infty}  ( \beta_n ' , \beta_n '') z^n$ be in
$\ThetaB$. For every $m \geq 1,$ consider the numbers $\gamma_m'$ and
$\gamma_m''$ defined as in the Equations (5.9) and (5.10) (in terms of
$\alpha_1', \alpha_1'', \ldots , \alpha_m', \alpha_m''$ and
$\beta_1', \beta_1'', \ldots , \beta_m', \beta_m'')$. Then the series
$\sum_{n=1}^{\infty} ( \gamma_n', \gamma_n'') z^n$ is called the boxed
convolution of type B of $f$ and $g,$ and is denoted
$f \ \freestarB \ g.$

\end{enumerate}

$\ $

The proposition proved above has as consequence that:

$\ $

\noindent
{\bf Corollary.} {\em The binary operation $\freestarB$ defined on
$\ThetaB$ is associative and has the series $\Delta ' (z) = (1,0)z$ as
a unit. }

\begin{proof} It is immediate that for every $n \geq 1$ it makes sense to
consider the truncation of $\freestarB$ to order $n$, thus obtaining a
binary operation $\freestarB_{\, n}$ on $( \C^2 )^n.$ (This is analogous to
the considerations done for $\freestarA$ in Section 5.1.) It clearly
suffices
to prove that, for every $n \geq 1,$ the operation $\freestarB_{\, n}$ is
associative and has the element
$( \ (1,0), (0,0), \ldots , (0,0) \ ) \in ( \C^2 )^n$ as a unit. But the
definition of $\freestarB$ is made so that we have
\begin{equation}
u_{ \alpha \freestarB_{\, n} \beta } \ \mid \ [ \ee , \omega ] \ = \
( \ u_{\alpha} \rstar u_{\beta} \ ) \ \mid \ [ \ee , \omega ],
\end{equation}
for every $n \geq 1$ and every $\alpha , \beta \in ( \C^2 )^n$. Hence, the
desired properties of $\freestarB_{\, n}$ follow from the corresponding ones
for the restricted convolution $\rstar$ (cf. Proposition 2.7).
\end{proof}

$\ $

%
%
%

\subsection
{A key connection between boxed convolutions of types A and B.}

In this section we present a way of relating the operations
$\freestarA$ and $\freestarB$, which will be crucial for understanding how
to move towards free probabilistic considerations of type B. In brief,
we will show that $\freestarB$ is still ``a $\freestarA$ operation,"
but with the scalars replaced by a certain algebra structure on $\C^2$.

$\ $

\noindent
{\bf Definition} \label{toeplitz algebra structure}
Let $\cal C$ denote $\C^{2}$ with the
multiplication given by
\begin{equation}
( \alpha ' , \alpha '') \cdot ( \beta ' , \beta '' ) \ := \
( \alpha ' \beta ' , \alpha ' \beta '' + \alpha '' \beta ' ).
\end{equation}

$\ $

Then ${\cal C}$ is a $\C$-algebra, with unit $(1,0)$.
This algebra may also be identified with the algebra of 2-by-2 upper
triangular Toeplitz matrices $\left[ \begin{array}{cl}
\alpha '  & \alpha '' \\
0         & \alpha '
\end{array} \right]$, or with $\C[x]/(x^2)$.
This algebra structure on
$\C^{2}$ is implicitly present in the considerations of
\cite{R}, Section 6, page 217.

Recall now that the boxed convolution of type A has a version
$\freestarA_{\, \cal C}$ which is ``with coefficients in ${\cal C}$''
(cf. the Remark and Notation at the end of Section 5.1). As a set, the
space of power series $\ThetaA_{\cal C}$ on which we consider the
operation $\freestarA_{\, \cal C}$ coincides with the set of power series
$\ThetaB$ of Definition 5.2. In other words, on $\ThetaB$ there are two
operations that we can look at: $\freestarB$, and $\freestarA_{\, \cal C}$.

$\ $

$\ $

\noindent
{\bf Theorem.} {\em In the notations of the preceding paragraph, we have
that $\freestarB = \freestarA_{\,\cal C}.$  }

\medskip

\begin{proof}
Let $f(z) = \sum_{n=1}^{\infty} ( \alpha_n ', \alpha_n '') z^n,$
$g(z) = \sum_{n=1}^{\infty} ( \beta_n ', \beta_n '') z^n$ be in $\ThetaB$.
Write
\[
f \ \freestarB \ g \ =: \
\sum_{n=1}^{\infty} ( \gamma_n ', \gamma_n '') z^n, \ \
\mbox{and} \ \ f \ \freestarA_{\, \cal C} \ g \ =: \
\sum_{n=1}^{\infty} ( \delta_n ', \delta_n '') z^n.
\]

We fix a positive integer $n$, for which we will show that
$( \gamma_n ' , \gamma_n '' )$ =
$( \delta_n ' , \delta_n '' ).$ It will be convenient to use the following
notations. First, for $p \in \nca (n)$ we will denote (as in
Section 1.1) the number of blocks of $p$ as $\blockno (p).$ Then given a
partition $p \in \nca (n),$ we will denote the blocks of $p$ (listed
in increasing order of their minimal elements, say) as
$F(p,i),$ $1 \leq i \leq \blockno (p)$. It will be moreover convenient
that, for every $p \in \nca (n),$ we use the notation
\[
F(p,i), \ \ \blockno (p) < i \leq n+1,
\]
for the blocks of the Kreweras complement $\Kr (p).$ (In other words
we set $F(p,i) :=$
\newline
$F( \Kr (p), i- \blockno (p) )$ for
$\blockno (p) < i \leq n+1.$ The indexing up to $n+1$ is correct by
virtue of the Eqn. (1.2) in Section 1.1.) Finally, it will also be
convenient to set a unified notation for the $\alpha$s and $\beta$s which
appear as coefficients for the series $f$ and $g$: for every $p \in \nca
(n)$
and $1 \leq i \leq n+1$ we put:
\[
\theta ' (p,i) \ := \ \left\{  \begin{array}{lll}
{\alpha '}_{\card( F(p,i) )}  & \mbox{if} & i \leq \blockno (p)  \\
{\beta '}_{\card( F(p,i) )}   & \mbox{if} & i >  \blockno (p),
\end{array}  \right.
\]
and
\[
\theta '' (p,i) \ := \ \left\{  \begin{array}{lll}
{\alpha ''}_{\card( F(p,i) )}  & \mbox{if} & i \leq \blockno (p)  \\
{\beta ''}_{\card( F(p,i) )}   & \mbox{if} & i >  \blockno (p).
\end{array}  \right.
\]

Let us look at $\gamma_n '$ and $\gamma_n ''$. First we have that:
\[
\gamma_n ' \ = \sum_{\begin{array}{c}
{\scriptstyle  p \in \nca (n) } \\
{\scriptstyle p := \{ F_1 , \ldots , F_k \}  }  \\
{\scriptstyle Kr(p) := \{ E_1 , \ldots , E_h \}  }
\end{array} }  \ 
\Bigl( \ \prod_{j=1}^{k} ( {\alpha '}_{\card(E_j)} )  \ \Bigr) \cdot
\Bigl( \ \prod_{i=1}^{h} {\beta '}_{\card(F_i)}  \ \Bigr)
\]
\begin{equation}
= \ \sum_{p \in \nca (n)} \ \prod_{i=1}^{n+1} \ {\theta '} (p,i).
\end{equation}
With $\gamma_n ''$ the situation would seem to be more complicated,
because this coefficient is defined via a summation over $\ncb (n)$,
as described by Equation (5.10) (where one replaces $m$ with $n$ in
Eqn.(5.10)). However, the summation over $\ncb (n)$ can be reduced
to one over $\nca (n)$, by using the $(n+1)$-to-1 cover
$\mbox{Abs} : \ncb (n) \rightarrow \nca (n)$ which is discussed in
Section 1.4. When doing so, and when taking into account the explicit
description of ${\mbox{Abs}}^{-1}(p)$ $(p \in \nca (n))$ provided by the
proof of Theorem 1.4, one gets:
\begin{equation}
\gamma_n '' \ = \ \sum_{p \in \nca (n)} \ \Bigl( \
\sum_{m=1}^{n+1} {\theta ''} (p,m) \cdot
\prod_{i \neq m} {\theta '} (p,i) \ \Bigr) .
\end{equation}

On the other hand, by recalling the definition of the operation
$\freestarA_{\, \cal C}$, we see that $( {\delta '}_n , {\delta ''}_n )$
equals:
\[ 
\sum_{\begin{array}{c}
{\scriptstyle  p \in \nca (n) } \\
{\scriptstyle p := \{ F_1 , \ldots , F_k \}  }  \\
{\scriptstyle Kr(p) := \{ E_1 , \ldots , E_h \}  }
\end{array} }  \ 
\Bigl( \ \prod_{j=1}^{k}
( {\alpha '}_{\card(E_j)} , {\alpha ''}_{\card(E_j)} )  \ \Bigr) \cdot
\Bigl( \ \prod_{i=1}^{h}
( {\beta '}_{\card(F_i)} , {\beta ''}_{\card(F_i)} )  \ \Bigr)
\]
\[
= \ \sum_{p \in \nca (n)} \ \prod_{i=1}^{n+1} \
( {\theta '} (p,i) , {\theta ''} (p,i) )
\]
(products considered with respect to the multiplication
on the algebra ${\cal C}$)
\begin{equation}
= \ \sum_{p \in \nca (n)} \ \Bigl( \
\prod_{i=1}^{n+1} \ {\theta '} (p,i), \
\sum_{m=1}^{n+1} {\theta ''} (p,m) \cdot
\prod_{i \neq m} {\theta '} (p,i) \ \Bigr) .
\end{equation}

By comparing (5.22) against (5.20)+(5.21), we obtain that
$( \gamma_n ' ,  \gamma_n '' )$ =
$( \delta_n ', \delta_n '' ),$ as desired.
\end{proof}

\medskip

\medskip

\medskip

\section
{ Non-crossing cumulants of type B.}

\medskip

In order to define the analogue of type B for non-crossing cumulants,
\setcounter{equation}{0}
we will pursue the idea of using the same equation as for
non-crossing cumulants of type A (Eqn.(4.2) in Section 4.3), but make
it have coefficients in $\C^{2}$, where $\C^{2}$ has the algebra structure
from Definition 5.3. We have to begin the discussion by introducing the
appropriate framework of non-commutative probability space.

\medskip

\medskip
\subsection
{Non-commutative probability space of type B.}

\medskip

{\bf Definition.} By a non-commutative probability space of type B
we will understand a system $\ncpsb$ where:

\begin{enumerate}
\item
$\A$ and $\varphi$ form a non-commutative probability space of type A.

\medskip
\item
$\V$ is a vector space over $\C$, and $f:  \V \rightarrow \C$ is a
linear functional. 

\medskip
\item
$\Phi : \A \times \V \times \A \rightarrow \V$ is a two-sided action
of $\A$ on $\V$. Usually we will simply write ``$a \xi b$'' instead of
$\Phi ( a , \xi , b)$, for $a,b \in \A$ and $\xi \in \V$.

\end{enumerate}

$\ $

\noindent
{\bf Remark.} Let $\ncpsb$ be a non-commutative probability space of
type B. On the vector space $\A \times \V$ we have a structure of unital
algebra (sometimes called ``the linking algebra of the bimodule $\V$''),
with multiplication defined by:
\begin{equation}
(a, \xi) \cdot (b, \eta ) \ := \ (ab, a \eta + \xi b), \ \
a,b \in \A, \ \xi , \eta \in \V .
\end{equation}
In other words, the algebra structure on $\A \times \V$ is the one obtained
when $( a, \xi ) \in \A \times \V$ is identified with a $2 \times 2$
matrix,
\[
(a , \xi ) \ \leftrightarrow \ \left[ \begin{array}{cc}
a & \xi \\ 0 & a
\end{array} \right] .
\]
The unit of $\A \times \V$ is $(I,0)$, where $I$ is the unit of $\A$.

Moreover, we have a natural linear map $E: \A \times \V \rightarrow
\C^{2}$,
defined by:
\begin{equation}
E( \ (a, \xi ) \ ) \ := \ ( \ \varphi (a), f( \xi ) \ ), \ \
a \in \A , \ \xi \in \V .
\end{equation}

\medskip

\medskip
\subsection{ Non-crossing cumulant functionals of type B.}

$\ $

\noindent
{\bf Definition.} Let $\ncpsb$ be a non-commutative probability space
of type B. The non-crossing cumulant functionals for this space are a
family of multilinear functionals $\Bigl( \kb_{n} :
( \A \times \V )^{n} \rightarrow \C^{2} \ \Bigr)_{n=1}^{\infty}$, uniquely
determined by the following equation: for every $n \geq 1$ and every
$a_{1}, \ldots , a_{n} \in \A$, $\xi_{1}, \ldots , \xi_{n} \in \V$
we have that:
\begin{equation}
\sum_{p \in \nca(n)} \ \prod_{ \begin{array}{c}
{\scriptstyle F \ block}  \\
{\scriptstyle of \ p}
\end{array}  } 
\kb_{\card (F)} \Bigl( \ (a_{1}, \xi_{1}), \ldots , (a_{n}, \xi_{n})
\mid F \ \Bigr)  =  E ( \ (a_{1}, \xi_{1}) \cdots (a_{n}, \xi_{n}) \ ).
\end{equation}
In (6.3) the following conventions of notation were used.

-- On the left-hand side of (6.3): if
$F = \{ j_{1} < j_{2} < \cdots < j_{m} \}$ is a subset of
$\{ 1, \ldots , n \}$, then
$( \ ( a_{1}, \xi_{1}), \ldots , (a_{n}, \xi_{n} ) \mid F \ )$ stands for
$( \ (a_{j_{1}}, \xi_{j_{1}}), \ldots , (a_{j_{m}}, \xi_{j_{m}}) \ )$
$\in \  ( \A \times \V )^{m}.$ The product indexed by ``$F$ block of $p$''
is considered with respect to the multiplication on $\C^{2}$
defined in Section 5.3.

-- On the right-hand side of (6.3), $E: \A \times \V \rightarrow \C^{2}$ is
as defined in Equation (6.2), and the product
$(a_{1}, \xi_{1}) \cdots (a_{n}, \xi_{n}) \in \A \times \V$ is computed
according to the multiplication rule of Equation (6.1). It is actually
straightforward to write the result of this multiplication explicitly:
\begin{equation}
(a_{1}, \xi_{1}) \cdots (a_{n}, \xi_{n})  \ = \
( a_{1} \cdots a_{n}, \sum_{m=1}^{n} a_{1} \cdots a_{m-1} \xi_{m}
a_{m+1}
\cdots a_{n} ).
\end{equation}

$\ $

$\ $

\noindent
{\bf Remarks.}

\vspace{6pt}

1. The Equation (6.3) defining the functionals $\kb_n$ has on its
left-hand side a summation over $\nca (n)$. However, when we concentrate
on the second component of this equation, we will really encounter a
summation over $\ncb (n),$ because of exactly the same phenomenon which
led to the formula $\freestarB = \freestarA_{\, \cal C}$ of Theorem 5.3.
This point will be re-appear in Section 6.5 below -- cf. the
remark at the end of that section.

\vspace{10pt}

2. As in type A, the recursive use of the Equation (6.3) gives explicit
formulas for the functionals $\kb_{n}$. These formulas repeat
those describing the functionals $\ka_{n}$ (as discussed in Section 4.3),
with the difference that now we operate in $\C^{2}$ instead of $\C$. For
instance:
\begin{equation}
\kb_{1} ( \ (a, \xi ) \ ) \ = \ E( \ (a, \xi ) \ ) \ = \
( \ \varphi (a), f( \xi ) \ ),
\end{equation}
while for $\kb_{2}$ the computation goes as follows:
\[
\kb_{2} ( \ (a_{1}, \xi_{1}), (a_{2}, \xi_{2}) \ ) \ = \
E( \ (a_{1}, \xi_{1}) (a_{2}, \xi_{2}) \ ) -
E( \ (a_{1}, \xi_{1}) \ ) E( \ (a_{2}, \xi_{2}) \ )
\]
\[
= \ E( \ (a_{1}a_{2}, a_{1} \xi_{2} + \xi_{1} a_{2} ) \ ) -
( \ \varphi (a_{1}) , f( \xi_{1} ) \ )
( \ \varphi (a_{2}) , f( \xi_{2} ) \ )
\]
\[
= \ ( \ \varphi (a_{1}a_{2}) , f(a_{1} \xi_{2}) + f(\xi_{1}a_{2}) \ ) -
( \ \varphi (a_{1}) \varphi (a_{2}) ,
\varphi (a_{1})f( \xi_{2} ) + f( \xi_{1} ) \varphi (a_{2}) \ )
\]
\begin{equation}
= \ \Bigl( \ \varphi (a_{1}a_{2}) - \varphi (a_{1}) \varphi (a_{2}) ,
\begin{array}[t]{r}
f(a_{1} \xi_{2}) - \varphi (a_{1})f( \xi_{2} ) \\
+ f(\xi_{1}a_{2}) \ ) - f( \xi_{1} ) \varphi (a_{2})
\end{array}
\ \Bigr).
\end{equation}

\vspace{10pt}

3. The equations defining the functionals $\kb_n$ are close to those used
in the framework of operator-valued cumulants developed in \cite{S2},
where 
$\C^2$ plays the role of algebra of scalars. There are however some
details which are different here, namely that

(a) $\C^2$ is not canonically embedded inside $\A \times \V ,$ and

(b) the map $E$ is not required to be a conditional expectation
(actually $\A \times \V$ doesn't even carry a canonical structure of
bimodule over $\C^2$).

\vspace{10pt}

4. Considering (6.6), it is natural to ask whether one couldn't also
describe the functionals $\kb_{n}$ by specifying  each component. The 
first component of $\kb_{2} ( \ (a_{1}, \xi_{1}), (a_{2}, \xi_{2}) \ )$ 
is just $\ka_{2} ( a_{1}, a_{2} )$. In order to describe what happens on 
the second component, we will have to discuss yet another variation of the
cumulant functional of type A, where one of the arguments is allowed to be
a vector. It is worth mentioning that exactly this variation of the
cumulant functional has recently appeared in \cite{NSS}, in connection
to the study of the concept of free Fisher information of Voiculescu.

\medskip

\medskip
\subsection
{The functionals $\kaa_{n;m}$.}

Let us observe that the Equation (4.2) defining the functional $\ka_{n}$
still makes sense when one  of the arguments
$a_{1}, \ldots , a_{n}$ of the functional is allowed to be a vector in an
$\A$--bimodule (rather than just being an element of $\A$).
We formalize this observation as follows.

$\ $

\noindent
{\bf Definition.} Let $\ncpsb$ be a non-commutative probability space of
type B. For every $n \geq 1$ and every $m \in \{ 1, \ldots , n \}$
we denote by $\kaa_{n;m}$ the multilinear functional from
$\A^{m-1} \times \V \times \A^{n-m}$ to $\C$ which is defined by exactly
the same formula as for $\ka_{n} : \A^{n} \rightarrow \C$, but where the
$m$th argument is a vector from $\V$, and where ``$\varphi$'' is replaced
by ``$f$'' in all the appropriate places.

$\ $

Referring to the explicit low order formulas presented in the Equations
(4.4) of Section 4.3, we have for instance that
$\kaa_{2;1} : \V \times \A \rightarrow \C$ and
$\kaa_{2;2} : \A \times \V \rightarrow \C$ are defined by:
\begin{equation}
\left\{  \begin{array}{lll}
\kaa_{2;1} ( \xi , a ) & = & f ( \xi a) - f ( \xi ) \varphi (a)  \\
                      &   &                                     \\
\kaa_{2;2} ( a , \xi ) & = & f ( a \xi ) - \varphi (a) f ( \xi ),
\end{array}  \right.
\end{equation}
for $\xi \in \V$ and $a \in \A$. Or that
$\kaa_{3;2} : \A \times \V \times \A \rightarrow \C$ is given by:
\[
\kaa_{3;2} ( a, \xi , a' ) \ = \ f(a \xi a') - \varphi (a) f( \xi a')
- f( \xi ) \varphi (aa') - f(a \xi ) \varphi (a') +
2 \varphi (a) f( \xi ) \varphi (a'),
\]
for $a,a' \in \A$ and $\xi \in \V$.
The general equation (analogous to (4.2) from Section 4.3) which determines
the functionals $\kaa_{n;m}$ can be written as follows. For $n \geq m \geq
1$
and for $a_{1}, \ldots , a_{m-1}, a_{m+1}, \ldots , a_{n}$
$\in \A$, $\xi \in \V$ we have:
\[
f( a_{1} \cdots a_{m-1} \xi a_{m+1} \cdots a_{n} ) \ =
\]
\begin{equation}
\sum_{p \in \nca (n)} \ \Bigl[ \  \kaa_{\card(F_{o});j}
( a_{1}, \ldots , a_{m-1}, \xi ,a_{m+1}, \ldots , a_{n} \mid F_{o} ) \cdot
\end{equation}
\[
\cdot \prod_{\begin{array}{c}
{\scriptstyle F \ block \ of \ p,}  \\
{\scriptstyle F \neq F_{o}}
\end{array}  } \ka_{\card(F)}
( a_{1}, \ldots , a_{m-1}, \xi ,a_{m+1}, \ldots , a_{n} \mid F ) \ \Bigr] ,
\]
where $F_{o}$ denotes the block of $p$ which contains $m$, and $j$
denotes the position of $m$ inside $F_{o}$ (i.e, if
$F_{o} = \{ i_{1} < i_{2} < \cdots < i_{l} \}$, then $i_{j}=m$).

$\ $

$\ $

\noindent
{\bf Remark} {\em (unified notation for $\ka_n$ and $\kaa_{m;n}$).}
Let $\ncpsb$ be a non-commutative probability space
of type B. In view of the fact that the cumulant functionals
$\ka_{n} : \A^{n} \rightarrow \C$ and
$\kaa_{n;m} : \A^{m-1} \times \V \times \A^{n-m} \rightarrow \C$ are
governed by exactly the same combinatorics, it is convenient
to use a unified notation for these two types of cumulants. More precisely,
it is convenient to deal with expressions of the form:
\begin{equation}
\mbox{``} \ \kaa_{n} ( x_{1}, \ldots , x_{n} ) \ \mbox{''},
\end{equation}
where:
\begin{equation}
\left\{ \begin{array}{cl}
\circ & \mbox{ Either all of $x_{1}, \ldots , x_{n}$ are in $\A$,}
\\
      & \mbox{ and $\kaa_{n} ( x_{1}, \ldots , x_{n} )  \ :=
            \ka_{n} (x_{1}, \ldots , x_{n}).$ }
\\
      &            
\\
\circ & \mbox{ Or there is a specified $m \in \{ 1, \ldots ,n \}$
              such that $x_{m} \in \V$, }
\\
      & \mbox{ while $x_1, \ldots , x_{m-1},x_{m+1}, \ldots ,x_n \in \A .$}
\\
      & \mbox{ In this case, $\kaa_{n} ( x_{1}, \ldots , x_{n} ) \ := \
                \kaa_{n;m} (x_{1}, \ldots , x_{n}).$ }
\end{array}  \right.
\end{equation}

$\ $

The use of the notation (6.9) simplifies for instance the above Equation
(6.8), in the sense that it is no longer necessary to distinguish the
block $F_o$ of $p$ in the expression on the right-hand side of (6.8).
Indeed, the term indexed by $p \in \nca (n)$ in the sum appearing there
is now simply written as:
\begin{equation}
\prod_{ F \ block \ of \ p} \kaa_{\card(F)}
(a_1, \ldots, a_{m-1} , \xi , a_{m+1}, \ldots , a_n \mid F ).
\end{equation}

Moreover, following a convention commonly used in the theory of
non-crossing cumulants of type A, it is convenient to introduce the
following notation.

$\ $

$\ $

\noindent
{\bf Notation.} Let $\ncpsb$ be a non-commutative probability space of
type B. Let $n$ be a positive integer, and let $p$ be a partition in
$\nca (n).$ For $x_1, \ldots , x_n$ as described in (6.10), we denote:
\begin{equation}
\kaa_p ( x_1 , \ldots , x_n ) \ := \
\prod_{ F \ block \ of \ p} \kaa_{\card(F)} (x_1, \ldots, x_n \mid F ).
\end{equation}

$\ $

$\ $

In the first situation covered by (6.10) (when $x_1, \ldots , x_n \in \A$),
the quantity appearing in (6.12) is precisely what is commonly denoted as
$\ka_p ( x_1, \ldots , x_n)$ in the theory of non-crossing cumulants of
type A. In the second situation covered by (6.10) (when $x_m \in \V$ for a
specified $m \in \{ 1, \ldots , n \}$), the quantity appearing in (6.12)
coincides with the one from (6.11).

Let us observe that if we use the above notations, then the Equation (6.8)
defining the cumulants of type $A'$ can be written in the more compressed
form
\begin{equation}
f( a_1 \cdots a_{m-1} \xi a_{m+1} \cdots a_n) \ = \
\sum_{p \in \nca (n)} \kaa_p (a_1, \ldots , a_{m-1}, \xi , a_{m+1},
\ldots , a_n ),
\end{equation}
for $n \geq m \geq 1$ and
$a_1, \ldots, a_{m-1}, a_{m+1}, \ldots , a_n \in \A$, $\xi \in \V$.

$\ $

We will conclude this subsection by recording two facts which are very
basic for the theory of cumulants of type A, and which have straightforward
generalizations to the cumulants of type $A'$.

$\ $

\noindent
{\bf Proposition.} {\em Let $\ncpsb$ be a non-commutative probability space
of type B. Let $n \geq 2$ be an integer, and consider some $x_1, \ldots
,x_n$
as described in (6.10). Then we have the relation: }
\[
\kaa_{n-1} (x_1, \ldots x_{r-1}, x_r x_{r+1}, x_{r+2}, \ldots , x_n)
= \ \kaa_{n} (x_1, \ldots , x_n)
\]
\begin{equation}
+ \ \sum_{ \begin{array}{c}
{\scriptstyle p \in \nca (n) \ with \ blno(p)=2,} \\
{\scriptstyle p \ separates \ r \ from \ r+1}
\end{array} } \ \ \kaa_p (x_1, \ldots , x_n).
\end{equation}  

$\ $

Note that on the left-hand side of (6.14) we may have that
$x_r x_{r+1} \in \A$ (if $x_r, x_{r+1} \in \A )$, or we may have that
$x_r x_{r+1} \in \V$ (if one of $x_r, x_{r+1}$ is from $\V$ while the other
one is from $\A$). The sum over $p$ appearing on the right-hand side of
(6.14) has $n-1$ terms, which could be easily listed explicitly
(the possible choices for $p$ are
$p= \{ (1, \ldots , r), \ (r+1, \ldots n) \} ,$
$p= \{ (l, \ldots , r), \ (1, \ldots , l-1, r+1, \ldots , n) \}$ for some
$2 \leq l \leq r,$ and
$p= \{ (1, \ldots , r, l+1, \ldots , n), \ ( r+1, \ldots , l) \}$
for some $r+1 \leq l \leq n-1$).

The relation (6.14) is copied from the theory of
cumulants of type A -- see \cite{S1}, or the generalizations obtained in
\cite{KS}. The proof of (6.14) is also copied verbatim from the type A
situation (cf. \cite{S1}, proof of Proposition 1 on page 622); this is due
to the fact that the cumulants $\kaa_{n}$ obey exactly the same
combinatorics as the cumulants $\ka_{n}$.

$\ $

$\ $

\noindent
{\bf Corollary.} {\em Let $\ncpsb$ be a non-commutative probability space of
type B. Let $n \geq 2$ be an integer, let $m$ be in $\{ 1, \ldots , n \}$,
and consider $\xi \in \V$ and
$a_1, \ldots ,a_{m-1}, a_{m+1}, \ldots , a_n \in \A .$
If there exists an index $r \in \{ 1, \ldots , n \} \setminus \{ m \}$
such that $a_{r} \in \C I,$ then
\begin{equation}
\kaa_{n} (a_1, \ldots , a_{m-1}, \xi , a_{m+1}, \ldots , a_n) = 0.
\end{equation} }

\medskip
\noindent
{\bf (Sketch of) Proof.} A well known fact in the theory of cumulants of
type A is that 
\begin{equation}
\ka_{n} (a_1, \ldots , a_n) = 0
\end{equation}
whenever $n \geq 2$ and there
exists $r \in \{ 1, \ldots ,n \}$ such that $a_r \in \C I$ (see e.g.
\cite{S1}, page 624). The proof of (6.15) is easily obtained by induction
on $n$, by using (6.16) and the recurrence relation (6.14). (Group $x_r$
with $x_{r+1}$ where $x_r = a_r \in \C I$. Then note that in the situation
at hand, most of the terms on the right-hand side of (6.14) will vanish,
either because of the induction hypothesis or because of (6.16).)
\qed

\medskip

\medskip

\subsection
{ Relation between the functionals $\kb_{n}$ and $\kaa_{n}$. }

We can now give the following alternative
description for the non-crossing cumulants of type B:

$\ $

\noindent
{\bf Theorem.} {\em Let $\ncpsb$ be a non-commutative probability space
of type B. Let $n$ be a positive integer, and consider the non-crossing
cumulant functional
$\kb_{n} : ( \A \times \V )^{n} \rightarrow \C^{2}$. Then we have:
\begin{equation}
\kb_{n} ( \ (a_{1}, \xi_{1}), \ldots , (a_{n}, \xi_{n} ) \ ) \ = \
\end{equation}
\[
( \ \ka_{n} (a_{1}, \ldots , a_{n} ) \ , \
\sum_{m=1}^{n} \kaa_{n} ( a_{1}, \ldots , a_{m-1}, \xi_{m},
a_{m+1}, \ldots , a_{n} ) \ ),
\]
for every $a_{1}, \ldots , a_{n} \in \A$ and
$\xi_{1}, \ldots , \xi_{n} \in \V$. }

\begin{proof} For every $n \geq 1$ let us denote by $\lambda_n$ the
multilinear functional from $( \A \times \V )^{n}$ to $\C^{2}$
defined by the right-hand side of (6.17). That is:
\begin{equation}
\lambda_{n} \Bigl( \ (a_{1} , \xi_{1}), \ldots ,
(a_{n} , \xi_{n}) \ \Bigr) \ := \
\end{equation}
\[
\Bigl( \ \ka_{n} ( a_1, \ldots , a_n ) , \
\sum_{m=1}^{n} \kaa_{n}
(a_1, \ldots a_{m-1}, \xi_{m}, a_{m+1}, \ldots a_{n} ) \ \Bigr) ,
\]
for $a_{1}, \ldots , a_{n} \in \A$, $\xi_{1}, \ldots , \xi_{n} \in \V$.
We want to show that $\lambda_{n} = \kb_{n}, \ \forall \ n \geq 1.$
In order to do so, it will suffice to verify that the Equation (6.3),
which determines uniquely the functionals $( \kb_n )_{n=1}^{\infty}$,
is also satisfied by the functionals $( \lambda_n )_{n=1}^{\infty}$.
In other words, it will suffice to show that:
\[
\sum_{p \in \nca (n)} \
\prod_{F \ block \ of \ p} \ \lambda_{\card(F)}
\Bigl( \ (a_1, \xi_1), \ldots , (a_n, \xi_n) \mid F \ \Bigr)
\]
\begin{equation}
= \ E \Bigl( \ (a_1, \xi_1), \ldots , (a_n, \xi_n) \ \Bigr),
\end{equation}
for every $n \geq 1$ and every $a_{1}, \ldots , a_{n} \in \A$,
$\xi_{1}, \ldots , \xi_{n} \in \V$.
We fix $n \geq 1$ and $a_{1}, \ldots , a_{n} \in \A$,
$\xi_{1}, \ldots , \xi_{n} \in \V$ about which we prove that (6.19) holds.

Let us also fix for the moment a partition $p \in \nca (n)$, and let
us focus on the term indexed by $p$ in the sum on the left-hand side
of (6.19). We compute:
\[
\prod_{ F \ block \ of \ p } \ \
\lambda_{\card(F)} ( \ (a_1, \xi_1), \ldots , (a_n, \xi_n ) \mid F \ )
\]
\[
= \ \prod_{ F \ block \ of \ p } \ \
\Bigl( \ \ka_{\card(F)} (a_1, \ldots , a_n \mid F), \
\sum_{m \in F} \kaa_{\card(F)}
(a_1, \ldots , a_{m-1}, \xi_m , a_{m+1}, \ldots , a_n \mid F ) \ \Bigr)
\]
(by substituting the two components of
$\lambda_{\card(F)} ( \ (a_1, \xi_1), \ldots , (a_n, \xi_n ) \mid F \ )$
from (6.18) )
\begin{equation}
= \ \Bigl( \ \prod_{\begin{array}{c}
{\scriptstyle F \ block} \\
{\scriptstyle of \ p}
\end{array} } \ \ka_{\card(F)} ( a_1 , \ldots , a_n \mid F), \ \
\sum_{m=1}^{n} \kaa_{p} (a_1, \ldots , a_{m-1}, \xi_m , a_{m+1}, \ldots ,
a_n) \ \Bigr)
\end{equation}
(by taking into account how the multiplication on $\C^2$ was defined in
Section 5.3, and by using the notation for $\kaa_p$ from Equation (6.12)
in Notation 6.3).

If we now sum over $p$ in (6.20), and if we take into account the
summation formulas described in Eqn.(4.2) of Section 4.3 (for the first
component) and in Eqn.(6.13) of Section 6.3 (for the second component),
then we find that the left-hand side of (6.19) is equal to:
\begin{equation}
\Bigl( \ \varphi ( a_1 \cdots a_n ), \ \sum_{m=1}^{n}
f( a_1 \cdots a_{m-1} \xi_{m} a_{m+1} \cdots a_n ) \ \Bigr) .
\end{equation}
But it is immediately seen that the right-hand side of (6.19) is also equal
to (6.21) (cf. Eqns.(6.4) and (6.2) in Section 6.1).
\end{proof}

\medskip

\medskip

\subsection{ Moment series and R-transform of type B.}

We can now introduce the type B analogues for the concepts of moment series
and of R-transform which were reviewed in the Section 4.4.

$\ $

\noindent
{\bf Definition.} Let $\ncpsb$ be a non-commutative probability space
of type B, and consider a couple $(a, \xi) \in \A \times \V$.
The moment series and R-transform of $(a, \xi )$ are the power series
$M$ and respectively $R$ defined as follows:
\begin{equation}
M(z) \ := \ \sum_{n=1}^{\infty} \ E( \ (a, \xi )^{n} \ ) z^{n};
\end{equation}
\begin{equation}
R(z) \ := \ \sum_{n=1}^{\infty} \
\kb_{n} ( \ \underbrace{(a, \xi ), \ldots , (a, \xi )}_{n} \ ) z^{n}.
\end{equation}
$M$ and $R$ belong to $\ThetaB$, the space of power series introduced in
Definition 5.2. The expressions ``$E( \ (a, \xi )^n \ )$'' appearing in
Equation (6.22) are computed according to the rules described in Remark 6.1.

$\ $

$\ $

A basic fact in the theory of the R-transform in type A is that the moment
series and the R-transform of a given element are related to each other
via a convolution formula (using the operation $\freestarA$). More
precisely: if $\ncps$ is a non-commutative probability space of type A, if
$a \in \A$, and if $M_a, R_a \in \ThetaA$ are the moment series and the
R-transform of $a$, then we have
\begin{equation}
M_a \ = \ R_a \ \freestarA \ \zeta ,
\end{equation}
where $\zeta (z) := \sum_{n=1}^{\infty} z^{n}$. The analogue of type B for
the Equation (6.24) is described as follows.

$\ $

\noindent
{\bf Proposition.} {\em Let $\ncpsb$ be a non-commutative probability
space
of type B, and consider a couple $(a, \xi) \in \A \times \V$. Then
the moment series $M$ and the R-transform $R$ of $(a, \xi )$ are related
by the formula:
\begin{equation}
M \ = \ R \ \freestarB \ \zeta '  ,
\end{equation}
where $\zeta ' \in \ThetaB$ is the series $\sum_{n=1}^{\infty} (1,0)z^{n}$.
}

\begin{proof} For every $n \geq 1$ we have that the coefficient of order
$n$ in $R \ \freestarB \ \zeta '$ is the same as the coefficient of order
$n$ in $R \ \freestarA_{\, \cal C} \ \zeta '$ (by Theorem 5.3). When
we express this coefficient by using the ${\cal C}$-valued version of
Equation (5.2), and when we take into account that all the coefficients
of $\zeta '$ are equal to the unit of ${\cal C}$, we obtain:
\[
\sum_{p \in \nca (n)} \ \prod_{F \ block \ of \ p}
\ \left( \begin{array}{c}
\mbox{coefficient of} \\
\mbox{order card$(F)$ in $R$}
\end{array} \right) .
\]
The latter quantity equals
\[
\sum_{p \in \nca (n)} \ \prod_{F \ block \ of \ p}
\ \kb_{\card(F)} \Bigl( \ (a, \xi ), \ldots , (a, \xi ) \ \Bigr)
\]
(by the definition of $R$ in Eqn.(6.23)), and is, hence, equal to
$E( \ (a, \xi )^{n} \ )$, by Equation (6.3). But this is
precisely the coefficient of order $n$ in $M$.
\end{proof}

$\ $

$\ $

\noindent
{\bf Remark.} The Equation (6.25) is equivalent to the particular case of
the formula (6.3) from Section 6.2, where the couples
$(a_1, \xi_1 ), \ldots , (a_n, \xi_n )$ (appearing in (6.3)) are all equal
to each other. If we had developed the more involved framework which
would allow us to define the operation $\freestarB$ for series in $k$
non-commuting indeterminants, then we could have now stated a more general
form of Equation (6.25), equivalent to the general case of (6.3).

The essential feature of (6.25) is that on the second component of
its right-hand side we encounter a summation over half of the lattice
$\ncb (n)$ -- or more precisely, over $\{ \pi \in \ncb (n) \ :$
$\pi$ has a zero-block $\}$. (The definition of $\freestarB$ would imply
a summation over all of $\ncb (n)$, but the terms indexed by $\pi$'s
without a zero-block are killed because of the special form of $\zeta '$.)
For instance taking the second component of the coefficient of $z^2$ on
the two sides of (6.25) yields the formula:
\[
f( a \xi + \xi a) \ = \
\Bigl[ \kaa_{2} ( \xi , a) + \kaa_{2} (a, \xi ) \Bigr] +
\Bigl[ \kaa_{1} ( \xi ) \ka_{1} (a) \Bigr] +
\Bigl[ \ka_{1} (a) \kaa_{1} ( \xi ) \Bigr] ,
\]
where the three groups of terms on the right-hand side are (respectively)
the contributions of the partitions $\{ \ ( 1,2,-1,-2) \ \} ,$
$\{ \ (1,-1), (2), (-2) \ \} ,$ $\{ \ (1),(-1), (2,-2) \ \} \in \ncb (2).$

\medskip

\medskip

\medskip

\section{ Free independence of type B.}

\medskip

Let $\ncpsb$ be a non-commutative probability space of type B.
\setcounter{equation}{0}
We will discuss the concept of free independence for a family
$( \A_{1}, \V_{1} ), \ldots , ( \A_{k}, \V_{k} )$, where
$\A_{1}, \ldots , \A_{k}$ are unital subalgebras of $\A$ and
$\V_{1}, \ldots , \V_{k}$ are linear subspaces of $\V$, such that
$\V_{j}$ is invariant under the (two-sided) action of $\A_{j}$,
$1 \leq j \leq k$. Since we have a type B analogue for non-crossing
cumulants,
our approach will go via the the counterpart of type B for the condition of
``vanishing of mixed cumulants'' (described in type A by the Equation (4.5)
of Section 4.3).

$\ $

\subsection{ Vanishing mixed cumulants in type B.}

\noindent
{\bf Definition.} Let $\ncpsb$ be a non-commutative probability space
of type B. Let $\A_{1}, \ldots , \A_{k}$ be unital subalgebras of $\A$
and let $\V_{1}, \ldots , \V_{k}$ be linear subspaces of $\V$, such that
$\V_{j}$ is invariant under the action of $\A_{j}$, for $1 \leq j \leq k$.
We say that $( \A_{1}, \V_{1} ), \ldots , ( \A_{k}, \V_{k} )$ have vanishing
mixed cumulants of type B if the following condition holds:
\begin{equation}
\left\{  \begin{array}{c}
\kb_{n} ( \ (a_{1}, \xi_{1}), \ldots , (a_{n}, \xi_{n}) \ ) = 0 \\
                                        \\
\mbox{whenever $a_{1} \in \A_{i_{1}}, \ldots , a_{n} \in \A_{i_{n}},
\xi_{1} \in \V_{i_{1}}, \ldots , \xi_{n} \in \V_{i_{n}}$ }             \\
                                        \\
\mbox{and $\exists \ 1 \leq s<t \leq n$ such that $i_{s} \neq i_{t}$.}
\end{array} \right.
\end{equation}

$\ $

$\ $

\noindent
{\bf Proposition.} {\em Let $\ncpsb$ and $\A_{1}, \ldots , \A_{k} \subset
\A$,
$\V_{1}, \ldots , \V_{k} \subset \V$ be as above, such that
$( \A_1, \V_1), \ldots , ( \A_k, \V_k)$ have vanishing cumulants of type B.
Then the following happen:

\medskip

1. $\A_{1}, \ldots , \A_{k}$ are freely independent in $\ncps$ (in the
type A sense).

\medskip

2. Suppose that $m,n \geq 0$ and that we have $m+n+1$ indices
$i_{m}, \ldots , i_{1},h,j_{1}, \ldots , j_{n}$ in $\{ 1, \ldots , k \}$
such that any two consecutive indices in this list are different from each
other ($i_{m} \neq i_{m-1}, \ldots ,i_1 \neq h \neq j_1, \ldots ,
j_{n-1} \neq j_{n}$). Suppose moreover that we have elements
$a_{m} \in \A_{i_{m}}, \ldots , a_{1} \in \A_{i_{1}}$,
$\xi \in \V_{h}$, 
$b_{1} \in \A_{j_{1}}, \ldots , b_{n} \in \A_{j_{n}}$, such that
\[
\varphi (a_{m}) = \cdots = \varphi (a_{1}) = 0 =  \varphi (b_{1}) =
\cdots = \varphi ( b_{n} ).
\]
Then:
\begin{equation}
f( a_{m} \cdots a_{1} \xi b_{1} \cdots b_{n} ) = 0 \
\mbox{ in the case when $m \neq n$,}
\end{equation}
and
\begin{equation}
f( a_{m} \cdots a_{1} \xi b_{1} \cdots b_{n} ) \ = \
\delta_{i_{1},j_{1}} \cdots \delta_{i_{n},j_{n}}
\varphi (a_{1}b_{1}) \cdots \varphi (a_{n}b_{n}) f( \xi )
\end{equation}
in the case when $m=n$. }

\begin{proof} 1. If we put $\xi_{1} = \cdots = \xi_{n} = 0$ in (7.1), and
look at the first component of the equality stated there, then we get
precisely the Equation (4.5) in Proposition 4.3. Hence, the free
independence
of $\A_{1}, \ldots , \A_{k}$ follows from Proposition 4.3.

\vspace{10pt}

2. We expand $f( a_m \cdots a_1 \xi b_1 \cdots b_n )$ as a sum, by
using Eqn.(6.13) from Section 6.3:
\[
f( a_m \cdots a_1 \xi b_1 \cdots b_n ) \ = \
\sum_{p \in \nca (m+n+1)} \ \kaa_p
(a_m, \ldots , a_1, \xi , b_1, \ldots , b_n )
\]
\begin{equation}
= \ \sum_{p \in \nca (m+n+1)} \ \Bigl( \ \prod_{F \ block \ of \ p} \
\kaa_{\card(F)} (a_m , \ldots , a_1, \xi , b_1 , \ldots b_n |F ) \
\Bigr) .
\end{equation}

Let us now pick a non-crossing partition $p \in \nca (m+n+1)$. We
consider conditions which force the the term indexed by $p$ in the sum on
the right-hand side of (7.4) to be  0.

\vspace{10pt} 

{\em Condition 1. $p$ has a singleton block $F = \{ s \}$ with $s \ne m+1$.}

Indeed, if $s < m+1$, then 
\[
\ka_{\card(F)} (a_m , \ldots , a_1, \xi , b_1 , \ldots b_n |F )
= \ka_{1} (a_s) = \varphi (a_s) = 0,
\]
and the term indexed by $p$ in (7.4) is 0. A similar argument applies
if $s > m+1.$

\vspace{10pt} 

{\em Condition 2. $p$ has a block $F$ which includes two successive
numbers.}

Indeed, since any two consecutive indices in the list
$(i_m, \dots, i_1, h, j_1, \dots, j_n)$ are different, the hypothesis of
vanishing of mixed cumulants of type B will give us in this case that
\begin{equation}
\kb_{m+n+1} ( \ (a_1,0), \ldots , (a_m,0), (0, \xi ), (b_1, 0), \ldots ,
(b_n, 0) \mid F \ ) \ = \ (0,0).
\end{equation}
If $F \not\ni m+1$ then we project the Equation (7.5) on its first component;
while if $F \ni m+1$, then we project the Equation (7.5) on its second
component. In either case, by using Theorem 6.4, we obtain that
\[
\kaa_{m+n+1} ( a_1, \ldots , a_m, \xi , b_1, \ldots , b_n \mid F ) \ = \ 0.
\]
This in turn implies that the term indexed by $p$ in (7.4) is equal to 0.

\vspace{10pt}

{\em Condition 3. $p$ has a block containing two distinct elements of
$X = \{1, \dots, m\}$ or of $Y = \{m+2, \dots, m + n + 1\}$. }

Indeed, in this case the non-crossing condition implies that $p$ has either
a singleton block or a block containing two successive numbers in $X$ or
$Y$.
Thus $p$ satisfies one of the conditions (1) or (2), and, therefore, the
summand indexed by $p$ in (7.4) is zero.

\vspace{10pt} 
 
Suppose $p$ is a non-crossing partition such that the term
indexed by $p$ in (7.4) is non-zero. Because of the non-crossing
condition, $p$ has a block which is a singleton or an interval of
length $\ge 2$. Because of the discussion concerning conditions (1) and
(2), such a block can only be the singleton $\{m+1\}$. It follows,
moreover, from the discussion concerning condition (3) that any other
block contains   exactly two elements, one from $\{ 1, \ldots , m \}$ and
the other from $\{ m+2 , \ldots , m+n+1 \}$.  But this is only possible
if $m = n$. Moreover, if $m = n$, the only non-crossing partition with
this block structure is $p = \{ (1,2n+1),(2,2n), \ldots, (n,n+2),(n+1) \}.$
  
Hence, in the case when $m \neq n,$ all the terms of the sum in (7.4)
vanish, and we obtain (7.2), while for $m=n,$ the sum on the right-hand
side of (7.4) reduces to only one term:
\begin{equation}
f( a_m \cdots a_1 \xi b_1 \cdots b_n ) \ = \
\ka_{2} (a_n , b_n ) \cdots \ka_{2} (a_1 , b_1 ) \cdot \kaa_{1} ( \xi ).
\end{equation}
It is immediate that
$\ka_{2} (a_r , b_r) = \delta_{i_{r},j_{r}} \varphi (a_r b_r),$ for every
$1 \leq r \leq n,$ and it is clear that
$\kaa_{1} ( \xi ) = f( \xi ),$ so that (7.3) follows from (7.6).
\end{proof}

\medskip

\subsection{ Free independence of type B, in terms of moments.}

The considerations of Section 7.1 prompt us to make the following
definition:

$\ $

\noindent
{\bf Definition.} Let $\ncpsb$ be a non-commutative probability space
of type B. Let $\A_{1}, \ldots , \A_{k}$ be unital subalgebras of $\A$
and let $\V_{1}, \ldots , \V_{k}$ be linear subspaces of $\V$, such that
$\V_{j}$ is invariant under the action of $\A_{j}$, for $1 \leq j \leq k$.
We will say that $( \A_{1}, \V_{1} ), \ldots , ( \A_{k}, \V_{k} )$ are 
freely independent if the following happen:

\medskip

(i) $\A_{1}, \ldots , \A_{k}$ are freely independent in $\ncps$ (in the
type A sense).

\medskip

(ii) We have the formula:
\begin{equation}
f( a_{m} \cdots a_{1} \xi b_{1} \cdots b_{n} ) \ = \
\left\{  
\begin{array}{ll}
0,                          & \mbox{ if } m \neq n \\
                           &                      \\
\delta_{i_{1},j_{1}} \cdots & \delta_{i_{n},j_{n}}
\varphi (a_{1}b_{1}) \cdots \varphi (a_{n}b_{n}) f( \xi ) \\
                           & \mbox{ if } m=n,
\end{array} \right.
\end{equation}
holding in the following context:

-- $m,n$ are non-negative integers;

-- $i_{m}, \ldots , i_{1},h,j_{1}, \ldots , j_{n}$ in $\{ 1, \ldots , k \}$
are such that any two consecutive indices in the list are different from
each
other; 

-- $a_{m} \in \A_{i_{m}}, \ldots , a_{1} \in \A_{i_{1}}$, $\xi \in \V_{h}$,
$b_{1} \in \A_{j_{1}}, \ldots , b_{n} \in \A_{j_{n}}$ are such that
$\varphi (a_{m}) = \cdots = \varphi (a_{1}) = 0$
$=  \varphi (b_{1}) = \cdots = \varphi ( b_{n} ).$

$\ $

$\ $

\noindent
{\bf Remark.} Let $\ncpsb$ be a non-commutative probability space of
type B, and let $E: \A \times \V \rightarrow \C^{2}$ be as defined in
Eqn.(6.2). Let $1 \in \A_1, \ldots , \A_k \subset \A$ be subalgebras, and
let $\V_1, \ldots , \V_k \subset \V$ be linear subspaces, such that $\V_{j}$ 
is invariant under the action of $\A_{j}$, for $1 \leq j \leq k.$ We denote:

$\A_o \ := \ \mbox{the subalgebra generated by }
\A_1 \cup \cdots \cup \A_k;$

$\V_o \ :=$ the smallest linear subspace of $\V$ which contains
$\V_1 \cup \cdots \cup \V_{k},$ and is invariant under the
action of $\A_{1}, \ldots , \A_{k}.$
\newline
If $(\A_1, \V_1), \ldots , (\A_k,\V_k)$ are freely independent,
then $E| \A_o \times \V_o$ is completely determined by the restrictions
$E| \A_1 \times \V_1, \ldots , E| \A_k \times \V_k.$

This statement amounts to two things.

(a) That $\varphi | \A_o$ is completely determined by
$\varphi | \A_1 , \ldots , \varphi | \A_k.$ This is a basic consequence of
the fact that $\A_1, \ldots , \A_k$ are freely independent in $\ncps$ in
the type A sense -- see \cite{VDN}, Section 2.5.

(b) That $f| \V_o$ is completely determined by
$\varphi | \A_1 , \ldots , \varphi | \A_k$ and by
$f | \V_1 , \ldots , f | \V_k.$ The phenomenon here is that $\V_o$ is the
linear span of vectors of the form $a_m \cdots a_1 \xi b_1 \cdots b_n,$
where  $a_m, \ldots , a_1, \xi, b_1, \ldots , b_n$ are as in part (ii) of
the preceding definition. (The straightforward verification of the latter
fact is left to the reader.)

$\ $

$\ $

\noindent
{\bf Proposition.} {\em Let $\ncpsb$ be a non-commutative probability space
of type B. Let $\A_{1}, \ldots , \A_{k}$ be unital subalgebras of $\A$ and
let $\V_{1}, \ldots , \V_{k}$ be linear subspaces of $\V$, such that $\V_{j}$
is invariant under the action of $\A_{j}$, for $1 \leq j \leq k$. If 
$( \A_1 , \V_1), \ldots , (\A_k, \V_k)$ are freely independent, then we have:
\begin{equation}
\left\{ \begin{array}{c}
\kaa_{n} ( a_1, \ldots , a_{m-1}, \xi , a_{m+1}, \ldots , a_n) = 0  \\
                      \\
\mbox{whenever $1 \leq m \leq n,$ $a_1 \in \A_{i_1}, \ldots ,
a_{m-1} \in \A_{i_{m-1}},$} \\
\mbox{  } \xi \in \V_{i_m}, a_{m+1} \in \A_{i_{m+1}},
\ldots , a_n \in \A_{i_n} \\
                      \\
\mbox{and $\exists \ 1 \leq s<t \leq n$ such that $i_s \neq i_t$.}
\end{array}   \right.
\end{equation} }

\begin{proof} By induction on $n.$ The case $n=1$ is fulfilled vacuously.
In the case $n=2$ we have to check that
\[
\kaa_{2} ( \xi , a) \ = \ \kaa_{2} ( a, \xi ) \ = \ 0
\]
when $a \in \A_i, \ \xi \in \V_j$ and $i \neq j$ $(1 \leq i,j \leq k).$
This follows immediately by using the hypothesis of free independence and
the concrete formulas for $\kaa_{2}( \xi , a ), \kaa_{2}(a , \xi )$ which
were mentioned in Section 6.3. Indeed, for instance for
$\kaa_{2} ( \xi , a ):$ the hypothesis of
free independence gives that $f( \ \xi ( a- \varphi (a) I ) \ ) \ = \ 0,$
hence that $f( \xi a ) \ = \ \varphi (a) f( \xi ),$ and then
\[
\kaa_{2} ( \xi , a) \
\stackrel{(6.7)}{=} \ f( \xi a) - f( \xi ) \varphi (a) \ = \ 0.
\]

The bulk of the proof will be devoted to the induction step: we fix
$n \geq 3,$ and we will prove that (7.8) is true for $n,$
by assumming that (7.8) was already proved for $1,2, \ldots , n-1.$

Consider some $m \in \{ 1, \ldots  n \}$ and some
$a_1 \in \A_{i_1}, \ldots , a_{m-1} \in \A_{i_{m-1}},$
$\xi \in \V_{i_m}, a_{m+1} \in \A_{i_{m+1}}, \ldots , a_n \in \A_{i_n},$
where $1 \leq i_1, \ldots , i_n \leq k$ and there exist
$1 \leq s<t \leq n$ such that $ i_s \neq i_t.$ We will treat separately
the following two cases:

\vspace{10pt}

{\em Case 1. There exists $r,$ $1 \leq r < n,$ such that $i_r = i_{r+1}.$ }

The Case 1 is treated by using the recurrence formula described in the
Equation (6.14) of Proposition 6.3, where we denote $x_{m} = \xi$ and
$x_r = a_r$ for $1 \leq r \leq n,$ $r \neq m,$ and where we group
$x_r$ with $x_{r+1}$ for a value of $r$ such that $i_r = i_{r+1}.$
Indeed, (6.14) will give us that:
\[
\kaa_{n} (a_1, \ldots , a_{m-1}, \xi , a_{m+1}, \ldots , a_n)
\ = \ \kaa_{n} ( x_1 , \ldots , x_n )
\]
\begin{equation}
= \ \kaa_{n-1} ( x_1, \ldots , x_{r-1}, x_r x_{r+1}, x_{r+2}, \ldots x_n )
- \ \sum_{ \begin{array}{c}
{\scriptstyle p \in \nca (n) \ with \ blno(p)=2 } \\
{\scriptstyle p \ separates \ r \ from \ r+1}
\end{array} }  \ \ \kaa_p ( x_1, \ldots , x_n ).
\end{equation}
Then each of the terms listed in (7.9) is found to be equal to 0, either
because of the induction hypothesis, or because of a phenomenon of vanishing
of mixed cumulants from $\A_1, \ldots , \A_k$ (in the type A sense).

\vspace{10pt}

{\em Case 2. We have that $i_r \neq i_{r+1},$ $\forall \ 1 \leq r<n.$ }

When treating the Case 2, note first that the value of
$\kaa_{n} (a_1, \ldots ,a_{m-1}, \xi , a_{m+1}, \ldots ,a_n )$ does
not change when we replace $a_r$ by $a_r - \varphi (a_r) I,$ for
$1 \leq r \leq n,$ $r \neq m;$ indeed, this is an
immediate consequence of the multilinearity of $\kaa_{n}$ combined
with the Corollary 6.3. By doing these replacements, we can assume without
loss of generality that $\varphi ( a_r ) = 0$ for all
$r \in \{ 1, \ldots ,n \} \setminus \{ m \} .$ But then the ``word''
$a_1 \cdots a_{m-1} \xi a_{m+1} \cdots a_n$ is exactly of the kind
considered in the Equation (7.7) from the definition of free independence,
and the hypothesis that $( \A_1 , \V_1 ), \ldots , ( \A_k , \V_k )$ are
freely independent gives us:
\begin{equation}
f( a_1 \cdots a_{m-1} \xi a_{m+1} \cdots a_n ) \ = \
\left\{ \begin{array}{ll}
0,                           &  \mbox{if $n \neq 2m-1$}  \\
                             &                           \\
\delta_{i_{1},i_{n}} \cdots  &  \delta_{i_{m-1},i_{m+1}} \varphi (a_1 a_n)
\cdots \varphi (a_{m-1}a_{m+1}) f( \xi )  \\
                             &   \mbox{if $n=2m-1$}
\end{array}  \right.
\end{equation}

On the other hand, we know from (6.13) that:
\begin{equation}
\kaa_{n} (a_1, \ldots , a_{m-1}, \xi , a_{m+1}, \ldots , a_n )
\ = \ f( a_{1} \cdots a_{m-1} \xi a_{m+1} \cdots a_{n} ) \ -
\end{equation}
\[
- \ \sum_{ \begin{array}{c}
{\scriptstyle p \in \nca (n) } \\
{\scriptstyle p \neq 1_{n}}
\end{array}  } \ \Bigl( \ \prod_{F \ block \ of \ p} \ \kaa_{\card(F)}
( a_{1}, \ldots , a_{m-1}, \xi ,a_{m+1}, \ldots , a_{n} \mid F ) \ \Bigr) .
\]
Most of the terms in the sum which was subtracted on the right-hand side of 
(7.11) are equal to 0. In fact, one can examine the three conditions listed
in the proof of  Proposition 7.1, and argue that, for $p$ satisfying at
least one of the three conditions, the term indexed by $p$ in the sum on the
right-hand side of (7.11) vanishes. Namely, if $p$ has a singleton block
$F = \{s\}$, where $s \ne m$, then
$$\ka_1(a_{1}, \ldots , a_{m-1}, \xi ,a_{m+1}, \ldots , a_{n} \mid F ) =
\ka_1(a_s) = \varphi(a_s) = 0.$$
In case $p$ has a block $F$ containing two successive numbers $s, s+1$, we
have
$$\kaa_{\card(F)}
( a_{1}, \ldots , a_{m-1}, \xi ,a_{m+1}, \ldots , a_{n} \mid F ) = 0, $$
by the induction assumption, since $\card(F) < n$ and $i_s \ne i_{s+1}$.
Finally, if $p$ has a block $F$ containing two elements of
$X =\{1, \dots, m-1 \}$ or of $Y = \{m+1, \dots, n\}$, then (exactly as
in the proof of Proposition 7.1) either $p$ has a singleton block
$F' = \{s\}$ with $s \ne m$, or $p$ has a block $F'$  containing two
successive numbers (and one of the two preceding arguments applies).
   
As observed in the proof of Proposition 7.1, the  non-crossing partitions
satisfying at least one of the three conditions will cover all non-crossing
partitions, if $n \neq 2m-1;$ and will cover all non-crossing partitions
with 
the exception of $p = \{ (1,n), \ (2,n-1), \ldots , (m-1,m+1), \ (m) \} ,$
if $n=2m-1.$ Thus from (7.11) we get:

\[
\kaa_{n} (a_1, \ldots , a_{m-1}, \xi , a_{m+1}, \ldots , a_n ) \ =
\]
\[
= \ \left\{  \begin{array}{l}
f(a_1 \cdots a_{m-1} \xi a_{m+1} \cdots a_n) , \mbox{ if } n \neq 2m-1 \\
                                          \\
f(a_1 \cdots a_{m-1} \xi a_{m+1} \cdots a_n)          \\
 - \ka_2 (a_1, a_n) \cdots \ka_2 (a_{m-1},a_{m+1}) \kaa_{1}( \xi ),
\mbox { if } n = 2m-1
\end{array}  \right.
\]
\vspace{10pt}
\[
= \ 0 \ \ \mbox{ (by Equation 7.10). }
\]
\end{proof}

$\ $

$\ $

\noindent
{\bf Corollary.} {\em Let $\ncpsb$ be a non-commutative probability space
of type B. Let $\A_{1}, \ldots ,$
$\A_{k}$ be unital subalgebras of $\A$ and let 
$\V_{1}, \ldots , \V_{k}$ be linear subspaces of $\V$, such that $\V_{j}$ 
is invariant under the action of $\A_{j}$, for $1 \leq j \leq k$. Then:
\begin{center}
$( \A_1, \V_1), \ldots , (\A_k, \V_k)$ are freely independent (in the
sense of Definition 7.2)

\vspace{10pt}

if and only if

\vspace{10pt}

$( \A_1, \V_1), \ldots , (\A_k, \V_k)$ have vanishing mixed cumulants
of type B 

(in the sense of Definition 7.1).
\end{center}
}

\begin{proof} ``$\Leftarrow$'' is the content of Proposition 7.1.

\vspace{10pt}

``$\Rightarrow$'' Let $(a_1, \xi_1) \in \A_{i_1} \times \V_{i_1},$
$\ldots , (a_n, \xi_n) \in \A_{i_n} \times \V_{i_n},$ with
$1 \leq i_1, \ldots , i_n \leq k,$ and suppose that there exist
$1 \leq s<t \leq n$ such that $i_s \neq i_t.$ We have to show that
$\kb_{n} ( \ (a_1, \xi_1), \ldots , (a_n, \xi_n) \ ) = 0.$ The two
components of $\kb_{n} ( \ (a_1, \xi_1), \ldots , (a_n, \xi_n) \ )$ are
$\ka_{n} ( a_1 , \ldots , a_n )$ and respectively $\sum_{m=1}^{n}
\kaa_{n} (a_1, \ldots , a_{m-1}, \xi_{m} , a_{m+1}, \ldots , a_n )$
(cf. Theorem 6.4). The first of the two components vanishes because of the
free independence of $\A_1, \ldots , \A_k$ in $\ncps$ (in the type A sense),
while the second of the two components vanishes by the preceding proposition.
\end{proof}

\medskip

\medskip

\subsection{ R-transforms for sums and products of free elements, in type B}

In order to make a case that (similarly to what we had in type A) the 
operation $\freestarB$ provides the middle-ground between the Cayley graph 
framework and free probability of type B, we now have to prove the analogue
of type B for the Theorem 5.1. We have:

$\ $

\noindent
{\bf Theorem.} {\em Let $\ncpsb$ be a non-commutative probability space of
type B. Let $\A_1 , \A_2$ be unital subalgebras of $\A$, let $\V_1 , \V_2$ be 
linear subspaces of $\V$, such that $\V_j$ is invariant under the action of 
$\A_j$ $(j = 1,2),$ and suppose that $( \A_1, \V_1)$, $( \A_2, \V_2)$ are
freely independent. Consider elements $a_1 \in \A_1, a_2 \in \A_2,$ 
$\xi_1 \in \V_1, \xi_2 \in \V_2,$ and denote the R-transform of 
$( a_j , \xi_j )$ by $R_j,$ for $j = 1,2.$ Then:
\begin{enumerate}
\item
The R-transform of $(a_1 , \xi_1) + (a_2 , \xi_2)$ is $R_1 + R_2.$
\item
The R-transform of $(a_1 , \xi_1) \cdot (a_2 , \xi_2)$ is 
$R_1 \ \freestarB \  R_2$ (where 
$(a_1 , \xi_1) \cdot (a_2 , \xi_2) := (a_1 a_2, a_1 \xi_2 + \xi_1 a_2),$ 
as discussed in Section 6.1).
\end{enumerate}  }

\begin{proof} Let us write explicitly 
$R_1 (z) = \sum_{n=1}^{\infty} v_n z^n$ and 
$R_2 (z) = \sum_{n=1}^{\infty} w_n z^n$, where
\[
v_n := \kb_n ( \ (a_1, \xi_1), \ldots , (a_1, \xi_1) \ ),  \ \
w_n := \kb_n ( \ (a_2, \xi_2), \ldots , (a_2, \xi_2) \ ),  \
\ \forall \ n \geq 1.
\]
The part (1) of the proposition amounts to the fact that 
\[
\kb_n ( \ (a_1, \xi_1) + (a_2, \xi_2), \ldots , (a_1, \xi_1) + 
(a_2, \xi_2) \ )  = v_n + w_n, \ \forall \ n \geq 1,
\]
which follows immediately from the multilinearity of the functionals 
$\kb_n$, and from the condition of vanishing mixed cumulants.

For the proof of part (2), let us denote the moment series and R-transform 
of $(a_1 , \xi_1) \cdot (a_2, \xi_2)$ as $M$ and $R$, respectively, and let
us denote the moment series of $(a_2, \xi_2)$ as $M_2.$ It is sufficient to
prove that 
\begin{equation}
M \ = \ R_1 \ \freestarB \ M_2.
\end{equation}
Indeed, in view of Proposition 6.5, the Equation (7.12) amounts to:
\begin{equation}
R \ \freestarB \ \zeta ' \ = \ R_1 \ \freestarB \ R_2 \ \freestarB \ \zeta ',
\end{equation}
and $\zeta '$ can be cancelled in (7.13) because it is invertible with respect 
to $\freestarB$.

By Theorem 5.3, we know that an equivalent form of (7.12) is 
\begin{equation}
M \ = \ R_1 \ \freestarA_{\, {\cal C} } \ M_2.
\end{equation}
On the left-hand side of (7.14), the coefficient of order $n$ is 
$E( \ ( \ (a_1, \xi_1) \cdot (a_2, \xi_2) \ )^n \ ).$ By writing
$( \ (a_1, \xi_1) \cdot (a_2, \xi_2) \ )^n$ as a product of $2n$ factors, and
by using the Equation (6.3) of Section 6.2, we write this coefficient as
\begin{equation}
\sum_{p \in \nca (2n)} \ \prod_{\begin{array}{c}
{\scriptstyle F \ block} \\
{\scriptstyle of \ p}  \end{array} } \ \ka_{\card (F)} \Bigl( \ 
(a_1, \xi_1), (a_2, \xi_2), \ldots , (a_1, \xi_1), (a_2, \xi_2) \mid F \ 
\Bigr) .
\end{equation}

Now, due to the condition of vanishing cumulants of type B for $(a_1, \xi_1)$ 
and $(a_2, \xi_2)$, we can in fact restrict the summation (7.15) to the 
partitions $p$ with the property that every block of $p$ is contained either
in $\{ 1,3, \ldots , 2n-1 \}$ or in $\{ 2,4, \ldots , 2n \} .$ A partition
$p \in \nca (2n)$ with the latter property is naturally parametrized by two
partitions $p_1, p_2 \in \nca (n)$ such that $p_2 \leq \Kr (p_1).$ (Namely,
$p$ is obtained by placing an isomorphic copy of $p_1$ on 
$\{ 1,3, \ldots , 2n-1 \}$, and an isomorphic copy of $p_2$ on 
$\{ 2,4, \ldots , 2n \} .$ The condition $p_2 \leq \Kr (p_1)$ is exactly the
one which ensures that $p$ is non-crossing as a partition of 
$\{ 1, \ldots , 2n \}$.) So in the end, the coefficient of order $n$ on the
left-hand side of (7.14) gets the form:
\[
\sum_{ \begin{array}{c}
{\scriptstyle p_1, p_2 \in \nca (n) } \\
{\scriptstyle p_2 \leq Kr (p_1) } 
\end{array} } \ \Bigl( \ \prod_{ \begin{array}{c}
{\scriptstyle F \ block} \\
{\scriptstyle of \ p_1} 
\end{array} } \ v_{\card (F)} \ \Bigr) \cdot
\Bigl( \ \prod_{ \begin{array}{c}
{\scriptstyle G \ block} \\
{\scriptstyle of \ p_2} 
\end{array} } \ w_{\card (G)} \ \Bigr) ,
\]
or in other words:
\begin{equation}
\sum_{p_1 \in \nca (n)} \ 
\Bigl( \ \prod_{ \begin{array}{c}
{\scriptstyle F \ block} \\
{\scriptstyle of \ p_1} 
\end{array} } \ v_{\card (F)} \ \Bigr) \cdot
\Bigr[ \ \sum_{ \begin{array}{c}
{\scriptstyle p_2 \in \nca (n) } \\
{\scriptstyle p_2 \leq Kr (p_1) } 
\end{array} } \ \prod_{ \begin{array}{c}
{\scriptstyle G \ block} \\
{\scriptstyle of \ p_2} 
\end{array} } \ w_{\card (G)} \ \Bigr] .
\end{equation}

Finally, by using again the Equation (6.3) from Section 6.2 (applied to the 
couple $(a_2, \xi_2)),$ and the definition of the operation 
$\freestarA_{\, {\cal C} }$, the quantity in (7.16) is easily identified as
the coefficient of order $n$ in $R_1 \ \freestarA_{\, {\cal C} } \ M_2.$
Thus the coefficients of order $n$ on the two sides of (7.14) are indeed 
equal to each other (for every $n \geq 1)$, and this completes the proof.
\end{proof}

$\ $

$\ $

Philippe Biane: Laboratoire de Math\'{e}matiques de l'Ecole Normale
Superieure, 45 rue d'Ulm, 75230 Paris, France.
Email: Philippe.Biane@ens.fr

$\ $

Frederick Goodman: Department of Mathematics, University of Iowa,
Iowa City, 
\newline
IA 52242, USA. 
Email: goodman@math.uiowa.edu.

$\ $

Alexandru Nica: Department of Pure Mathematics, University of Waterloo,
\newline
Waterloo, Ontario N2L 3G1, Canada.
Email: anica@math.uwaterloo.ca


\medskip

\medskip

\medskip

\begin{thebibliography}{99}

\bibitem{B1} P. Biane.
Minimal factorizations of a cycle and central multiplicative functions
on the infinite symmetric group, Journal of Combinatorial Theory
Series A 76 (1996), 197-212.

\bibitem{B2} P. Biane.
Some properties of crossings and partitions, Discrete Mathematics
175 (1997), 41-53.

\bibitem{BS} M. Bozejko, R. Speicher.
$\psi$-independent and symmetrized white noises, in Quantum Probability
and Related Topics (L. Accardi editor), World Scientific, Singapore,
VI (1991), 219-236.

\bibitem{DRS} P. Doubilet, G.-C. Rota, R. Stanley.
On the foundations of combinatorial theory (VI): The idea of
generating function, Proceedings of the sixth Berkeley symposium on
mathematical statistics and probability, Lucien M. Le Cam et al.
editors, University of California Press, 1972, 267-318.

\bibitem{dlH} P. de la Harpe.
Topics in geometric group theory, University of Chicago Press, 2000.

\bibitem{H} J.E. Humphreys.
Reflection groups and Coxeter groups, Cambridge University Press, 1990.

\bibitem{KS}
B. Krawczyk, R. Speicher.
Combinatorics of free cumulants, Journal of Combinatorial Theory Series A
90 (2000), 267-292.

\bibitem{K} G. Kreweras.
Sur les partitions non-crois\'{e}es d'un cycle, Discrete Mathematics
1 (1972), 333-350.

\bibitem{NS1} A. Nica, R. Speicher.
A ``Fourier transform'' for multiplicative functions on non-crossing
partitions, Journal of Algebraic Combinatorics 6 (1997), 141-160.

\bibitem{NS2} A. Nica, R. Speicher.
On the multiplication of free $n$-tuples of non-commutative random
variables. With an Appendix by D. Voiculescu: Alternative proofs
for the type II free Poisson variables and for the compression results.
American Journal of Mathematics 118 (1996), 799-837.

\bibitem{NS3} A. Nica, R. Speicher.
The combinatorics of free probability, Preliminary version,
December 1999. Lecture notes of Centre \'{E}mile Borel, the
Henri Poincar\'{e} Institute, Paris, France).

\bibitem{NSS} A. Nica, D. Shlyakhtenko, R. Speicher.
Operator-valued distributions I. Characterizations of freeness,
to appear in International Mathematics Research Notices.

\bibitem{R} V. Reiner.
Non-crossing partitions for classical reflection groups,
Discrete Mathematics 177 (1997), 195-222.

\bibitem{S1} R. Speicher.
Multiplicative functions on the lattice of non-crossing partitions
and free convolution, Mathematische Annalen 298(1994), 611-628.

\bibitem{S2} R. Speicher.
Combinatorial theory of the free product with amalgamation and
operator-valued free probability theory, Memoirs of the Amer. Math. Soc.
132 (1998), x+88.

\bibitem{V1} D. Voiculescu.
Addition of certain non-commutative random variables, Journal of
Functional Analysis 1986.

\bibitem{VDN} D. Voiculescu, K. Dykema, A. Nica.
Free random variables, CRM Monograph Series, volume 1, AMS, 1992.

\end{thebibliography}
\end{document}